\documentclass[reqno]{amsart}
\usepackage{amsmath}
\usepackage{amssymb}
\usepackage{amscd}
\usepackage{graphics}
\usepackage{latexsym}
\usepackage{stmaryrd}
\usepackage{hyperref}

\theoremstyle{plain}
\newtheorem{theorem}{Theorem}[section]
\newtheorem{thm}[theorem]{Theorem}
\newtheorem{cor}[theorem]{Corollary}
\newtheorem{lem}[theorem]{Lemma}
\newtheorem{prop}[theorem]{Proposition}
\newtheorem{claim}[theorem]{Claim}

\theoremstyle{definition}
\newtheorem{defn}[theorem]{Definition}

\newtheorem{ques}[theorem]{Question}

\newtheorem{rmk}[theorem]{Remark}

\newtheorem{notat}[theorem]{Notation}

\newtheorem{hyp}[theorem]{Hypothesis}

\theoremstyle{remark}

\newcommand{\ZZ}{\mathbb{Z}}

\newcommand{\PP}{\mathbb{P}}

\newcommand{\mc}{\mathcal}

\newcommand{\OO}{\mc{O}}

\newcommand{\marpar}[1]{}

\newcommand{\M}{\overline{\mc{M}}}

\newcommand{\lt}{\left}
\newcommand{\rt}{\right}

\newsavebox{\sembox}
\newlength{\semwidth}
\newlength{\boxwidth}

\newsavebox{\semrbox}
\newlength{\semrwidth}
\newlength{\boxrwidth}

\title[Very twisting lines on Grassmannians]
{Very twisting families of pointed lines on Grassmannians}

\author[de Jong]{A. J. de Jong}
\address{Department of Mathematics \\ 
Massachusetts Institute of Technology \\ 
Cambridge MA 02139}
\email{dejong@math.mit.edu} 

\author[Starr]{Jason Michael Starr}
\address{Department of Mathematics \\
Massachusetts Institute of Technology \\
Cambridge MA 02139}
\email{jstarr@math.mit.edu} 

\date{\today}

\begin{document}


\begin{abstract}
This excerpt is a section from an article in progress.  This section
proves that for the \emph{Grassmannians}, i.e., the homogeneous spaces
of Picard number one for the classical simple algebraic groups, there
exists a very twisting family of pointed lines.
\end{abstract}


\maketitle



\section{Very twisting lines on Grassmannians} \label{sec-grasslines}
\marpar{sec-grasslines}

\medskip\noindent
This is an extract from an article in progress relating rational
connectedness of spaces of rational curves to existence of sections of
families over surfaces.  Two other articles, ~\cite{dJS7} and
~\cite{dJS8}, also deal with aspects of this work.

This extract is concerned with a very limited
problem: extending the basic arguments from ~\cite{HS2} to
Grassmannians and isotropic Grassmannians.  This might
seem superfluous since Kim and Pandharipande prove rationality of 
the spaces of
rational curves on every projective
homogeneous spaces, ~\cite{KP}.  However, for sections of families
over surfaces, one needs also the existence of very twisting families
of lines, which is what this note proves. 

\medskip\noindent
Let $\kappa$ be an algebraically closed field.  
Let $(X,\OO_X(1))$ be a quasi-projective $\kappa$-variety together
with an  
ample invertible sheaf.  Denote by $X_\text{sm}$ the smooth locus of
$X$.  The scheme $\M_{0,1}(X_\text{sm},1)$ represents the functor of
pointed lines in $X$.  On $\M_{0,1}(X_{\text{sm}},1)$, there is a rank
$2$ locally free sheaf $E$, a rank $1$ locally direct summand $L$ of
$E$, and a morphism $g:\PP(E)\rightarrow X_\text{sm}$ pulling back
$\OO_{X_\text{sm}}(-1)$ to the universal rank $1$ locally direct
summand $\OO_{\PP(E)}(-1)$ of $\pi^*E$.  The universal line is the
$\PP^1$-bundle $\pi:\PP(E)\rightarrow \M_{0,1}(X_\text{sm},1)$, the
universal section of $\pi$, $\sigma:\M_{0,1}(X_\text{sm},1)\rightarrow
\PP(E)$, pulls back $\OO_{\PP(E)}(-1)$ to the locally direct summand
$L$, and the universal map is $g$.

\medskip\noindent
On $\M_{0,1}(X_{\text{sm}},1)$ there is an important invertible sheaf
$\psi^\vee$, defined as $\sigma^*\OO_{\PP(E)}(\text{Image}(\sigma)) =
\sigma^* T_\pi$, the pullback of the normal bundle of $\sigma$, which
is also the pullback of the vertical tangent bundle of $\pi$.  
Equivalently, $\psi^\vee$ satisfies a canonical isomorphism,
$$
\psi^\vee \cong \text{det}(E)\otimes (L^\vee)^{\otimes 2}.
$$

\medskip\noindent
The \emph{evaluation morphism}, denoted
$\text{ev}:\M_{0,1}(X_\text{sm},1)\rightarrow X_\text{sm}$, is $g\circ
\sigma$.  The open subset where $\text{ev}$ is smooth is denoted $U$.  On
this open set, an important locally free sheaf is the vertical tangent
bundle $T_\text{ev}$ of $\text{ev}$, i.e., the dual of the sheaf of
relative differentials of $\text{ev}$.  
The morphism $(\pi,g):\PP(E) \rightarrow \M_{0,1}(X_\text{sm},1)\times
X_\text{sm}$ is a regular embedding, and thus has a locally free
normal sheaf $N$.  
The open subset $U$ is the maximal open set
over which $N$ is
$\pi$-relatively globally generated, and $T_\text{ev}$ has an
equivalent definition,
$$
T_\text{ev} \cong \pi_*(N(-\text{Image}(\sigma))).
$$  

\begin{defn}\cite{HS2} \label{defn-vt}
\marpar{defn-vt}
A \emph{very twisting family of pointed lines on $X$} is a morphism
$\zeta:\PP^1 \rightarrow \M_{0,1}(X_{\text{sm}},1)$ such that,
\begin{enumerate}
\item[(i)]
$U$ contains $\text{Image}(\zeta)$,
\item[(ii)]
$\zeta^* T_{\text{ev}}$ is ample, and
\item[(iii)]
the degree of $\zeta^*\psi^\vee$ is nonnegative.
\end{enumerate}
\end{defn}

\medskip\noindent
Using the canonical isomorphisms, 
(i)--(iii) are equivalent to,
\begin{enumerate}
\item[(i')]
the restriction of $N$ to the fiber of $\pi$ over each point of
$\text{Image}(\zeta)$ is globally generated,
\item[(ii')]
$\zeta^*\pi_*(N(-\sigma(\PP^1)))$ is ample, and
\item[(iii')]
$\zeta^*\text{det}(E)\otimes(L^\vee)^{\otimes 2}$ has nonnegative degree.
\end{enumerate}

\medskip\noindent
The open subset $U$ intersects
each irreducible component of
$\M_{0,1}(X_\text{sm},1)$
whose lines cover a dense subset of $X_\text{sm}$, 
cf. 
\cite[1.1]{KMM92c}. 

\begin{lem} \label{lem-hom}
\marpar{lem-hom}
If $T_{X_\text{sm}}$ is globally generated, then $U$ equals
$\M_{0,1}(X_{\text{sm}},1)$.  In particular, if $X$ is a homogeneous
space $G/P$, then $U$ equals $\M_{0,1}(X,1)$.
\end{lem}

\begin{proof}
Since $T_{X_\text{sm}}$ is globally generated, $g^*T_{X_\text{sm}}$ is
globally generated, and thus the quotient $N$ is globally generated.
For a homogeneous space $G/P$, $T_X$ is globally generated by $T_e
G\otimes_\kappa \OO_X$.
\end{proof}

\medskip\noindent
Unfortunately, there typically exist rational
curves in $\M_{0,1}(G/P,1)$ on which
$\psi^\vee$ or $T_{\text{ev}}$ has negative degree.
Thus, if there exists a very twisting family, it is a special rational
curve.  
There are
some obvious special rational curves in $\M_{0,1}(X,1)$; in some cases
these give very twisting families.

\begin{defn} \label{defn-orbit}
\marpar{defn-orbit}
Let $\lambda:\mathbf{G}_m\rightarrow G$ be a 1-parameter subgroup.
Let $s:G\times \M_{0,1}(G/P,1)\rightarrow \M_{0,1}(G/P,1)$ be the
canonical action.  Let $p\in \M_{0,1}(G/P,1)$ be a point.  There is an
induced morphism $\zeta^\text{o}:\mathbf{G}_m \rightarrow \M_{0,1}(G/P,1)$ by 
$\zeta(t) = s(\lambda(t),p)$.  
This extends uniquely to a morphism $\zeta:\PP^1
\rightarrow \M_{0,1}(G/P,1)$, by the valuative criterion of properness.
An \emph{orbit curve} is a morphism
$\zeta$ thus obtained.
\end{defn}

\begin{ques} \label{ques-grasslines}
\marpar{ques-grasslines}
Does there exist a very twisting family
$\zeta:\PP^1\rightarrow\M_{0,1}(G/P,1)$?  
\end{ques}

\medskip\noindent
If the answer is affirmative, a second question is whether there
exists a very twisting orbit curve.

\medskip\noindent
We answer Question~\ref{ques-grasslines} when $X$ is the Grassmannian
$\text{Flag}(k,V)$ of rank $k$ subspaces of an $n$-dimensional vector
space $V$ and 
when $X$ is the Grassmannian $\text{Flag}_k(V,\beta)$ of rank $k$
isotropic subspaces of an $n$-dimensional vector space $V$ with a
symmetric or 
skew-symmetric bilinear pairing $\beta$.
Thus we answer the question 
when $G$ is one of the
classical simple groups $\mathbf{SL}_n$, $\mathbf{SO}_n$,
$\mathbf{Sp}_{2n}$ and $P$ is a maximal parabolic group.

\medskip\noindent
\textbf{The exceptional cases.}
There are some exceptional cases: there does not exist a very twisting
family of pointed lines if $X$ equals a finite set, $\PP^1$, or
$\PP^1\times \PP^1$.  
In these case $\text{ev}$ is finite, and thus
$T_\text{ev}$ is the zero sheaf.  
\begin{enumerate}
\item[(i)] 
For the classical Grassmannian, the single exceptional case is
$(n,k)= (2,1)$.
\item[(ii)]
In the skew-symmetric case,
the single exceptional case is $(n,k) = (2,1)$.
\item[(iii)]
In the symmetric case, the exceptional cases are
$(n,k) = (2,1),(3,1), (4,1)$, and $(4,2)$.  
\end{enumerate}

\begin{thm} \label{thm-main}
\marpar{thm-main}
For every pair of positive integers $(n,k)$ satisfying $n\geq 2k$ and
not on the exceptional list above,
there is a very twisting family of pointed lines to $X$. 
In many cases, there is a very twisting orbit curve.
\end{thm}

\section{Pointed lines on classical Grassmannians} \label{subsec-pl}
\marpar{subsec-pl}
Let $V$ be a rank $n$ $\kappa$-vector space, and let $k$ be an
integer, $0<k<n$.  Denote by $X$ the Grassmannian $\text{Flag}(k,V)$,
and denote by $S(k,V)$ the universal rank $k$ locally direct summand
of $V\otimes_\kappa \OO_X$.  Denote by $\OO_X(1)$ the ample invertible
sheaf giving the Pl\"ucker embedding, i.e., the ample generator of the
Picard group of $X$.  Denote by $\M$ the scheme $\M_{0,1}(X,1)$.  
Denote by $M$ the flag variety $\text{Flag}(k-1,k,k+1;V)$, and denote
by $S_{k-1}\subset S_k\subset S_{k+1}\subset V\otimes\OO_M$ the
universal $(k-1,k,k+1)$-flag of locally direct summands of $V$.

\begin{lem} \label{lem-gg}
\marpar{lem-gg}
On $\PP(E)$, the locally free sheaf $g^*S(k,V)^\vee$ is globally
generated.  On $\M$, $\pi_*[g^*S(k,V)^\vee\otimes \OO_{\PP(E)}(-1)]$ is
an invertible sheaf, and $\pi_*[g^*S(k,V)^\vee]$ is locally free of
rank $k+1$.  The tautological map $V^\vee\otimes\OO_{\M}\rightarrow
\pi_*[g^*S(k,V)^\vee]$ is surjective.
\end{lem}

\begin{proof}
Since $V^\vee\otimes\OO_X$ generates
$S(k,V)^\vee$, $V^\vee\otimes \OO_{\PP(E)}$ generates
$g^*S(k,V)^\vee$.  So the restriction of $g^*S(k,V)^\vee$ to every
fiber of $\pi$ is isomorphic to $\OO_{\PP^1}(a_1)\oplus\dots \oplus
\OO_{\PP^1}(a_k)$ for integers $0\leq a_1\leq \dots \leq a_k$.   
By definition, $\OO_X(1) =
\bigwedge^k S(k,V)^\vee$ has degree $1$ on every fiber of $\pi$.  Thus
the restriction of $g^* S(k,V)^\vee$ to every fiber of $\pi$ is
isomorphic to $F:= \OO_{\PP^1}(1)\oplus \OO_{\PP^1}^{\oplus(k-1)}$.  
So, firstly, $\pi_*[g^*S(k,V)^\vee]$ is locally free of rank
$h^0(\PP^1,F)=k+1$, and $\pi^*[g^*S(k,V)^\vee\otimes
  \OO_{\PP(E)}(-1)]$ is locally free of rank $h^0(\PP^1,F(-1)) = 1$.  
Since the
only subspace of $H^0(\PP^1,F)$ generating $F$ is all of
$H^0(\PP^1,F)$, and since $V^\vee\otimes \OO_{\PP(E)}$ generates
$g^*S(k,V)^\vee$, also $V^\vee\otimes \OO_{\M} \rightarrow
\pi_*[g^*S(k,V)^\vee]$ is surjective.
\end{proof}

\medskip\noindent
There is a $(1,2)$-flag of locally direct summands of
$\pi_*[g^*S(k,V)^\vee]$,
$$
\pi_*[g^*S(k,V)^\vee\otimes \OO_{\PP(E)}(-1)]\otimes (E/L)^\vee
\subset \pi_*[g^*S(k,V)^\vee \otimes \OO_{\PP(E)}(-1)]\otimes E^\vee
\subset \pi_*[g^*S(k,V)^\vee].
$$
Dually, there is a $(k-1,k)$-flag of locally direct summands of
$(\pi_*[g^*S(k,V)^\vee])^\vee$,
$$
\begin{array}{r}
(\pi_*[g^*S(k,V)^\vee\otimes \OO_{\PP(E)}(-1)]\otimes E^\vee)^\perp
\subset \\
(\pi_*[g^*S(k,V)^\vee \otimes \OO_{\PP(E)}(-1)]\otimes
(E/L)^\vee)^\perp \subset \\
(\pi_*[g^*S(k,V)^\vee])^\vee.
\end{array}
$$
Because $V^\vee\otimes \OO_M\rightarrow \pi_*[g^*S(k,V)^\vee]$ is
surjective, $(\pi_*[g^*S(k,V)^\vee])^\vee$ is canonically a rank $k+1$
locally
direct summand of $V\otimes\OO_{\M}$.  This defines a $(k-1,k,k+1)$-flag
of locally direct summands of $V\otimes\OO_{\M}$, denoted
$E_{k-1}\subset E_k \subset E_{k+1}\subset V\otimes \OO_{\M}$,
$$
\begin{array}{ccc}
E_{k-1} & := &
(\pi_*[g^*S(k,V)^\vee\otimes \OO_{\PP(E)}(-1)]\otimes E^\vee)^\perp, \\
E_k     & := &
(\pi_*[g^*S(k,V)^\vee \otimes \OO_{\PP(E)}(-1)]\otimes
(E/L)^\vee)^\perp, \\ 
E_{k+1} & := &
(\pi_*[g^*S(k,V)^\vee])^\vee.
\end{array}
$$
By the
universal property of the flag variety, there exists a unique morphism
$\iota':\M \rightarrow M$ pulling back
$S_{k-1}\subset S_k \subset S_{k+1}\subset V\otimes\OO_M$ to 
$E_{k-1}\subset E_k\subset E_{k+1}\subset
V\otimes\OO_{\M}$.  

\begin{prop} \label{prop-iota}
\marpar{prop-iota}
The morphism $\iota'$ is an isomorphism.  Moreover,
$\text{ev}:\M\rightarrow \text{Flag}(k,V)$ is the composition of
$\iota'$ with the tautological projection
$\text{Flag}(k-1,k,k+1;V)\rightarrow \text{Flag}(k,V)$.
\end{prop}

\begin{proof}
On $M$, denote by
$E'$ the rank $2$ locally free sheaf $S_{k+1}/S_{k-1}$.  Denote by
$L'$ the rank $1$ locally direct summand $S_k/S_{k-1}$.  Denote by
$\pi':\PP(E')\rightarrow \text{Flag}(k-1,k,k+1;V)$ the associated
$\PP^1$-bundle.  
There is a
unique section $\sigma':\text{Flag}(k-1,k,k+1;V)\rightarrow \PP(E')$
pulling back the universal rank $1$ locally direct summand
$\OO_{\PP(E')}(-1)$ of $(\pi')^*E'$ to $L'$.  
On $\PP(E')$,
there is a rank $k$ locally direct summand $S'_k$ of $V\otimes
\OO_{\PP(E')}$ defined to be the preimage in $(\pi')^*S_{k+1}$ of the
universal rank $1$ locally direct summand $\OO_{\PP(E')}(-1)\subset
(\pi')^*(S_{k+1}/S_{k-1})$.  
By the universal property of $X$, there
exists a unique morphism $g':\PP(E')\rightarrow X$ pulling back the
locally direct summand $S(k,V)$ to $S'_k$.  
By definition of $\sigma'$,
$(\sigma')^*S'_k$ equals $S_k$ as a locally direct summand of
$V\otimes\OO_F$.  Therefore $g'\circ \sigma'$ is the tautological
projection $\text{Flag}(k-1,k,k+1;V)\rightarrow \text{Flag}(k,V)$.

\medskip\noindent
By the definition of $\OO_X(1)$, $(g')^*\OO_X(-1)$ is isomorphic to
$\bigwedge^k S'_k$.  By definition of $S'_k$, this is isomorphic to
$(\pi')^*\bigwedge^{k-1}S_{k-1}\otimes \OO_{\PP(E')}(-1)$.  In
particular, the restriction of $(g')^*\OO_X(1)$ to every fiber of
$\pi'$ is isomorphic to $\OO_{\PP^1}(1)$.  Thus $(\pi',\sigma',g')$ is
a family of pointed lines in $X$.    
By the universal property of $\M_{0,1}(X,1)$, there exists a unique
morphism $\iota:M\rightarrow \M$ pulling back
$(\pi,\sigma,g)$ to $(\pi',\sigma',g')$.  It follows easily that
$\iota'\circ \iota$ is the identity map $\text{Id}_M$.  

\medskip\noindent
To prove that
$\iota\circ \iota'$ is $\text{Id}_{\M}$, it suffices to find
an isomorphism $h:\PP(E)\rightarrow (\iota')^*\PP(E')$ such that
$(\iota')^*\sigma'$ equals $h\circ\sigma$ and $g$ equals $(\iota')^*g'\circ
h$.  By construction of $\iota'$, there is a canonical isomorphism of
$(\iota')^*E'$ with $E\otimes
(\pi_*[g^*S(k,V)^\vee\otimes \OO_{\PP(E)}(-1)])^\vee$.  
By the universal property of $\PP(E')$, there exists a unique
isomorphism $h:\PP(E)\rightarrow 
(\iota')^*\PP(E')$ pulling back the universal locally direct summand
$(\iota')^*\OO_{\PP(E')}(-1)$ of $(\iota')^*(\pi')^*E'$ to the locally
direct summand,
$$
\OO_{\PP(E)}(-1) \otimes
(\pi^*\pi_*[g^*S(k,V)\otimes\OO_{\PP(E)}(-1)])^\vee 
\subset
\pi^* E\otimes 
(\pi^*\pi_*[g^*S(k,V)\otimes\OO_{\PP(E)}(-1)])^\vee \cong \pi^*(\iota')^*E',
$$
obtained from the universal locally direct summand
$\OO_{\PP(E)}(-1)\subset \pi^* E$.  
Since
$(\iota')^*L'$ equals
$L\otimes (\pi_*[g^*S(k,V)^\vee\otimes
  \OO_{\PP(E)}(-1)])^\vee$ as locally direct summands of $(\iota')^*E'$,
$(\iota')^*\sigma'$ equals $h\circ \sigma$.  To prove $(\iota')^*g' \circ h$
equals $g$, it suffices to prove that $(\iota')^*S'_k$ equals $g^*S(k,V)$ as
locally direct summands of $V\otimes \OO_{\PP(E)}$.  By definition of
$\iota'$, the locally direct summand $(\iota')^*(\pi')^*S_{k-1}$ of
$V\otimes\OO_{\PP(E)}$ equals $\pi^*E_{k-1}$, which is contained in
$g^*S(k,V)$.  Forming the corresponding quotients, 
it suffices to prove that $g^*S(k,V)/\pi^*E_{k-1}$ equals
$(\iota')^* \OO_{\PP(E')}(-1)$ as locally direct summands of
$(\iota')^*(\pi')^*E'$, i.e., of $\pi^*(E_k/E_{k-1})$.  
As $g^* S(k,V)/\pi^* E_{k-1} \subset \pi^*(E_k/E_{k-1})$ is
isomorphic to,
$$
\OO_{\PP(E)}(-1) \otimes
(\pi^*\pi_*[g^*S(k,V)\otimes\OO_{\PP(E)}(-1)])^\vee 
\subset
\pi^* E\otimes 
(\pi^*\pi_*[g^*S(k,V)\otimes\OO_{\PP(E)}(-1)])^\vee,
$$
compatibly with the isomorphism to $\pi^*\iota^*E'$, this follows from
the definition of $h$.
\end{proof}

\begin{cor} \label{cor-iota}
\marpar{cor-iota}
Denote by $\iota:\text{Flag}(k-1,k,k+1;V)\rightarrow M$ the inverse
morphism of $\iota'$.  There are canonical isomorphisms,
$$
\iota^*T_{\text{ev}}\cong [(S_{k+1}/S_k)^\vee\otimes ((V\otimes
\OO_F)/S_{k+1})] \oplus [(S_k/S_{k-1})\otimes S_{k-1}^\vee],
$$ 
$$
\iota^*\psi^\vee \cong (S_{k+1}/S_k)\otimes (S_k/S_{k-1})^\vee.
$$
\end{cor}

\begin{proof}
Because $\text{ev}\circ \iota'$ is the tautological projection,
$(\iota')^*T_{\text{ev}}$ is the vertical tangent bundle of the
projection.  This projection is the fiber product of the relative
Grassmannian $\text{Flag}(1,(V\otimes \OO_F)/S_k)$ and the
relative Grassmannian $\text{Flag}(k-1,S_k)$.  
The first isomorphism follows from the
well-known computation of the vertical tangent bundle of a
Grassmannian bundle.  The second isomorphism follows from the
isomorphisms $(\iota')^*L = L' = S_k/S_{k-1}$ and $(\iota')^*(E/L) =
E'/L' = S_{k+1}/S_k$.  
\end{proof}

\section{Pointed lines on isotropic Grassmannians} \label{subsec-BCD1}
\marpar{subsec-BCD1}
Assume that $\text{char}(\kappa)$ is not $2$.
This subsection gives the analogues of Proposition~\ref{prop-iota} and
Corollary ~\ref{cor-iota} in the isotropic case.  
Let $V$ be an $n$-dimensional vector space with a symmetric or
skew-symmetric nondegenerate pairing.  Let $X$ be the Grassmannian of
isotropic $k$-planes in $V$.
The later sections prove existence of a very twisting family of
pointed lines.  The proof breaks up into
several cases.   
\begin{enumerate}
\item[I.]
This is the case when $k$ is odd and $n\geq \max(4,2k+2)$.
\item[II.]
This is the case when $k$ is even and $n\geq 2k+2$.  
\item[III.]
This is the case when $n=2k$ and the pairing is symmetric.  The case
when $n=2k-1$ and the pairing is symmetric reduces to this case.
\item[IV.]
This is the case when $n=2k$ and the pairing is skew-symmetric.
\end{enumerate}

\begin{defn} \label{defn-pairing}
\marpar{defn-pairing}
Let $B$ be a $\kappa$-scheme.  A \emph{symmetric pairing over $B$}, resp. a
\emph{skew-symmetric pairing on $B$}, is a triple $(E,\mc{L},\beta)$ of a
locally  free $\OO_B$-module $E$ of finite, constant rank, an
invertible $\OO_B$-module $\mc{L}$, and an
isomorphism of $\OO_B$-modules,
$\beta:E\rightarrow
E^\vee\otimes \mc{L}$ 
such that
$\beta$ equals $\beta^\dagger\otimes \text{Id}_{\mc{L}}$,
resp. $\beta$ equals
$-\beta^\dagger\otimes \text{Id}_{\mc{L}}$.
If the invertible sheaf $\mc{L}$ equals $\OO_B$, the pair is written
$(E,\beta)$.  
\end{defn}

\begin{notat} \label{notat-BCD}
\marpar{notat-BCD}
Let $(E,\beta)$ be a symmetric pairing over $B$ or a skew-symmetric
pairing over $B$.  
For every increasing sequence of
integers $\underline{k} = (k_1<\dots<k_r)$, denote by
$\text{Flag}_{\underline{k}}(E,\beta)$ the bundle over $B$ parametrizing
$\underline{k}$-flags of isotropic locally direct summands of $E$.  Denote by
$\pi:\text{Flag}_{\underline{k}}(E,\beta)\rightarrow B$ the projection,
and denote by $S_{k_1}(E,\beta)\subset \dots \subset
S_{k_r}(E,\beta)\subset \pi^*E$ the universal $\underline{k}$-flag of
isotropic locally direct summands. 
\end{notat} 

\medskip\noindent
Let $k$ and $n$ be positive integers with $n\geq 2k$.
Let $W$ be a $\kappa$-vector space of dimension $n$.  Let $(W,\beta)$
be a symmetric or skew-symmetric pairing.
Denote by $X$ the isotropic flag variety,
$$
X = \text{Flag}_k(W,\beta).
$$
By the universal property of $\text{Flag}(k;W)$, there exists a unique
morphism $e:X\rightarrow \text{Flag}(k;W)$ pulling back the universal
locally direct summand $S(k,W)$ of $W$ to the universal isotropic
locally direct summand $S_k(W,\beta)$.  The morphism $e$ is a closed
immersion.  Except when $\beta$ is odd and $n=2k$ or $2k-1$, the
invertible sheaf $\OO_X(1)$ is defined to be the pullback by $e$ of
the Pl\"ucker invertible sheaf $\OO(1)$ on $\text{Flag}(k;W)$.  In
these cases, the induced morphism
$\M_{0,1}(e,1):\M_{0,1}(X,1)\rightarrow \M_{0,1}(\text{Flag}(k;W),1)$
is a closed immersion.  Thus, using Proposition~\ref{prop-iota}, there
is a canonical closed immersion of $\M_{0,1}(X,1)$ in
$\text{Flag}(k-1,k,k+1;W)$.  The cases when $n=2k$ or $2k-1$ are a bit
more complicated.  

\medskip\noindent
\textbf{Symmetric case, $n\geq 2k+2$.}  
Denote by $\M$ the scheme $\M_{0,1}(X,1)$.
Denote by $M$ the isotropic flag variety,
$$
M = \text{Flag}_{k-1,k,k+1}(W,\beta).
$$
Denote by $S_{k-1}\subset S_k\subset S_{k+1} \subset W\otimes_\kappa
\OO_M$ the universal isotropic flag.  By the universal property of
$\text{Flag}(k-1,k,k+1;W)$, there is a unique morphism
$e':M\rightarrow \text{Flag}(k-1,k,k+1;W)$ pulling back the universal
flag to $S_{k-1}\subset S_k\subset S_{k+1}\subset W\otimes_\kappa
\OO_M$.  The morphism $e'$ is a closed immersion.

\medskip\noindent
Define $E'$ to be the rank $2$
locally free sheaf, $S_{k+1}/S_{k-1}$, and define $L'$ to be the
invertible sheaf, $S_k/S_{k-1}$.  Denote by $\pi':\PP(E')\rightarrow M$
the associated $\PP^1$-bundle.  There is a unique section
$\sigma':M\rightarrow \PP(E')$ such that the pullback of the rank $1$
locally direct summand $\OO_{\PP(E')}(-1)$ of $(\pi')^*E'$ equals
$L'$.  On $\PP(E')$ there is a rank $k$ locally  direct
summand $S'_k$ of $(\pi')^* S_{k+1}$ 
defined as the preimage of the rank $1$
locally direct summand $\OO_{\PP(E')}(-1)$ of
$(\pi')^*(S_{k+1}/S_{k-1})$.  Because $(\pi')^*S_{k+1}$ is a rank $k+1$
locally direct summand of $W\otimes_\kappa \OO_{\PP(E')}$, $S'_k$ is a
rank $k$ locally direct summand of $W\otimes_\kappa \OO_{\PP(E')}$.
Because 
$(\pi')^* S_{k+1}$ is isotropic, also $S'_k$ is isotropic.  By the
universal property of $\text{Flag}_k(W,\beta)$, there exists a unique
morphism $g':\PP(E')\rightarrow \text{Flag}_k(W,\beta)$ such that $(g')^*
S_k$ equals $S'_k$ as a locally direct summand of $W\otimes_\kappa
\OO_{\PP(E')}$.  

\begin{prop} \label{prop-flagBD}
\marpar{prop-flagBD}
Assume $n$ is at least $2k+2$.
The datum $(\pi',\sigma',g')$ is a family of pointed
lines in $\text{Flag}_k(W,\beta)$ parametrized by $M$.  The associated
morphism $\iota:M\rightarrow \M$ is an
isomorphism, compatible with the closed immersions into
$\text{Flag}(k-1,k,k+1;W)$.   Moreover, $\text{ev}\circ \iota$ equals
the tautological projection
$\text{Flag}_{k-1,k,k+1}(W,\beta)\rightarrow \text{Flag}_k(W,\beta)$.
\end{prop}

\begin{proof}
The proof that $(\pi',\sigma',g')$ is a family of pointed lines is
identical to the argument in the proof of
Proposition~\ref{prop-iota}.  The morphism $\iota$ is clearly
compatible with the closed immersions into
$\text{Flag}(k-1,k,k+1;W)$.  So $\iota$ is a closed immersion.  Since
both $M$ and $\M$ are smooth, to prove $\iota$ is an isomorphism, it
suffices to prove that $\iota$ is surjective.  The fact about
$\text{ev}$ follows from the analogous fact in Proposition~\ref{prop-iota}.

\medskip\noindent
Let $[E_{k-1}\subset
E_k\subset E_{k+1}\subset W]$ be any flag in $\text{Flag}(k-1,k,k+1;W)$
contained in the image of $\M$.  Because the pointed line in
$\text{Flag}(k;W)$ associated to this flag is contained in
$\text{Flag}_k(W,\beta)$, 
for every rank $k$ subspace
$S'_k\subset E_{k+1}$ containing $E_{k-1}$, $S'_k$ is isotropic.
In particular, since every vector in $E_{k+1}$ is contained in such a
subspace, every vector in $E_{k+1}$ is isotropic.  By the polarization
identity for symmetric bilinear pairings (which holds because
$\text{char}(k)$ is not two!), the subspace $E_{k+1}$ is isotropic.
Thus $[E_{k-1}\subset E_k\subset E_{k+1}\subset W]$ is contained in
the image of $M$.  
\end{proof}

\medskip\noindent
On $M$ there is a rank $n-k-1$ locally direct summand $S_{k+1}^\perp$
of $W\otimes_\kappa \OO_M$ defined as the annihilator of $S_{k+1}$
under $\beta$.  Since $S_{k+1}$ is isotropic, by definition $S_{k+1}$
is a locally direct summand of $S_{k+1}^\perp$. 

\begin{cor} \label{cor-flagBD}
\marpar{cor-flagBD/prop-flagBD2}
Assume $n$ is at least $2k+2$.
The pullbacks under $\iota$ of $T_\text{ev}$ and $\psi^\vee$ admit
canonical isomorphisms,
$$
\iota^*T_{\text{ev}}\cong [(S_{k+1}/S_k)^\vee\otimes
(S_{k+1}^\perp/S_{k+1})] \oplus [(S_k/S_{k-1})\otimes S_{k-1}^\vee], 
$$ 
$$
\iota^*\psi^\vee \cong (S_{k+1}/S_k)\otimes (S_k/S_{k-1})^\vee.
$$
\end{cor}

\begin{proof}
The projection $\text{Flag}_{k-1,k,k+1}(W,\beta)\rightarrow
\text{Flag}_k(W,\beta)$ is the fiber product of 
the relative isotropic Grassmannian,
$\text{Flag}_1(S_k^\perp/S_k,\widetilde{\beta})$, and      
the relative classical Grassmannian, $\text{Flag}(k-1,S_k)$.  The
proof of the corollary is almost identical the the proof of
Corollary~\ref{cor-iota}.  The one new element is
the well-known isomorphism of the vertical tangent bundle of
$\text{Flag}_1(E,\beta)\rightarrow B$ with $S_1(E,\beta)^\vee\otimes
(S_1(E,\beta)^\perp/S_1(E,\beta))$, using Notation~\ref{notat-BCD}.
\end{proof}

\begin{rmk} \label{rmk-flagBD3}
\marpar{rmk-flagBD3}
If $k=1$, then $S_{k-1}$ is the zero sheaf.  If $n=2k+2$, then
$S_{k+1}^\perp$ equals $S_{k+1}$ so that $S_{k+1}^\perp/S_{k+1}$ is
the zero sheaf.
\end{rmk}

\medskip\noindent
\textbf{Symmetric case, $n=2k+1,2k+2$.}  If $n$ equals $2k+1$ or
$2k+2$, then the Picard group is $\ZZ$, but the invertible sheaf
giving the Pl\"ucker embedding, the \emph{Pl\"ucker invertible sheaf}, 
is not a generator of the Picard group.  First consider the
case when $n=2k+1$.  Let $(\mathbf{1},\beta_{\mathbf{1}})$ be a
symmetric pairing such that $\text{dim}(\mathbf{1})$ equals $1$.   
Define $(W',\beta')$ to be the orthogonal direct sum of  $(W,\beta)$
and
$(\mathbf{1},\beta_{\mathbf{1}})$.  This has rank $n'=n+1$.  Denote
$k+1$ by $k'$. 
Denote by $M'$ the isotropic flag variety,
$M' = \text{Flag}_{k'}(W',\beta')$.  Because every isotropic subspace
of $W$ has dimension $\leq k$, no $k'$-dimensional isotropic
subspace of $W'$ is contained in $W$.  Therefore the following
morphism of $\OO_{M'}$-modules is surjective,
$$
S_{k'} \rightarrow W'\otimes_\kappa \OO_{M'}
\xrightarrow{\text{pr}_{\mathbf{1}}}
\mathbf{1}\otimes_\kappa \OO_{M'}.
$$
Denote the kernel by $K$.  This is a rank $k$ locally direct
summand of $W\otimes_\kappa \OO_{M'}$.  Because $S_{k'}$ is isotropic
for $\beta'$, $K$ is isotropic for $\beta$.  
By the universal property of $\text{Flag}_k(W,\beta)$, there exists a
unique morphism
$e:M' \rightarrow M$,
pulling back the universal isotropic flag $S_k \subset W\otimes_\kappa
\OO_M$ to $K\subset W\otimes_\kappa \OO_{M'}$.  

\begin{lem} \label{lem-flagB4}
\marpar{lem-flagB4}
If $n=2k+1$, the morphism $e:\text{Flag}_{k+1}(W',\beta')\rightarrow
\text{Flag}_k(W,\beta)$ is an \'etale, finite morphism of degree $2$ 
identifying $\text{Flag}_{k+1}(W',\beta')$ with a disjoint union of 2
copies of $\text{Flag}_k(W,\beta)$.
\end{lem}

\begin{proof}
The first part 
is just the fact that a symmetric, bilinear pairing on a rank $2$
vector space has precisely $2$ isotropic lines.  That this cover is
trivial can be checked directly.  It also follows from the fact that
$\text{Flag}_k(W,\beta)$ is separably rationally connected, in fact
separably unirational, together with a corollary of Koll\'ar:  a
separably rationally 
connected variety has trivial \'etale fundamental group,
cf. ~\cite[Cor. 3.6]{DBX}.
\end{proof}

\medskip\noindent
Of course $\text{dim}(W') = n+1 = 2(k+1) = 2k'$.  Thus the case $n=2k+1$ is 
reduced to the case $n'=2k'$.  

\medskip\noindent
Next consider the case when $n=2k$.  Let $X$ be one of the two
connected components of $\text{Flag}_k(V,\beta)$.  As above, there is
an embedding of $X$ in $\text{Flag}(k,V)$, and thus a Pl\"ucker
invertible sheaf on $X$.  Because $X$ is smooth and rational, there is
no torsion in the Picard group of $X$.  Therefore there exists a
minimal ample invertible sheaf some power of which equals
the Pl\"ucker invertible sheaf.

\begin{notat} \label{notat-gen}
\marpar{notat-gen}
Assume $n$ equals $2k$.
Denote by $\OO_X(1)$ the unique minimal ample invertible sheaf on
$\text{Flag}_k(W,\beta)$ some power of which equals the Pl\"ucker
invertible sheaf.
\end{notat}

\begin{lem} \label{lem-O1}
\marpar{lem-O1}
The Picard group of $X$ is generated by $\OO_X(1)$ and the Pl\"ucker
invertible sheaf is isomorphic to $\OO_X(2)$.  In particular, a
smooth rational curve
in $X$ is a line with respect to $\OO_X(1)$ iff it has degree $2$ with
respect to the Pl\"ucker invertible sheaf.
\end{lem}

\begin{proof}
Write $X$ as $\textbf{SO}_{2k}/P = \textbf{Spin}_{2k}/P'$.  There is a
natural isomorphism of the Picard group of $X$ and the character group
of $P'$.  Choose an isomorphism $W \cong V\oplus V^\vee$ sending
$\beta$ to the canonical symmetric bilinear pairing on $V\oplus
V^\vee$.  Then the stabilizer group $P$ of the isotropic 
flag $V\subset V\oplus V^\vee$ is the group of all maps,
$$
\lt(
\begin{array}{c|c}
U^{-1} & B \\
\hline
0 & U^\dagger 
\end{array} \rt),
$$
where $U$ is in $\textbf{GL}(V)$ and $B:V^\vee \rightarrow V$ is any
skew-symmetric map.  The group $P'$ is a connected extension of $P$ by
$\boldmath{\mu}_2$.  Because $\text{char}(\kappa)$ is not $2$, there is
no nontrivial $\boldmath{\mu}_2$-extension of the additive group of
skew-symmetric matrices.  Therefore $P'$ is the basechange by
$P\rightarrow \textbf{GL}(V)$ of a connected
$\boldmath{\mu}_2$-extension $G'$ of $\textbf{GL}(V)$.  There is precisely
one such, namely the basechange by
$\text{det}:\textbf{GL}(V)\rightarrow \textbf{G}_m$ of the unique
connected $\mathbf{\mu}_2$ extension of $\mathbf{G}_m$,
$(*)^2:\mathbf{G}_m\rightarrow \mathbf{G}_m$.  The character group of
$P'$ equals the character group of $G'$, which is a free Abelian group
containing the character group of $\textbf{GL}(V)$ as an index $2$
subgroup.  The Pl\"ucker invertible sheaf corresponds to the character
$\text{det}:\textbf{GL}(V)\rightarrow \textbf{G}_m$.  Since this
character is twice the generator of the character group of $G'$, the
Pl\"ucker invertible sheaf is isomorphic to the square of the
generator of the Picard group of $X$. 
\end{proof}

\medskip\noindent
Denote by $\M$ the the space $\M_{0,1}(X,1)$.  
Denote by $M$ the isotropic flag variety,
$$
M = \text{Flag}_{k-2,k}(W,\beta).
$$
Denote by $S_{k-2}\subset S_k \subset W\otimes_\kappa \OO_M$ the
universal isotropic flag.  Denote by $S_{k-2}^\perp$ the rank $k+2$
locally direct summand of $W\otimes_\kappa \OO_M$ defined as the
annihilator of $S_{k-2}$.  Because $S_k$, resp. $S_{k-2}$, is
isotropic, it is a 
locally direct summand of $S_{k-2}^\perp$.  
Denote by $F$ the rank $4$ locally free sheaf, $S_{k-2}^\perp/S_{k-2}$.
Associated to $\beta$ there is a symmetric pairing $(F,\beta_F)$ on
$M$.
Denote by $C'$ the isotropic flag variety,
$$
C' = \text{Flag}_2(E,\widetilde{\beta}).
$$
Denote by $\pi:C'\rightarrow M$ the projection.  
Denote by $G$ the rank $2$, locally direct summand 
$S_k/S_{k-2}$ of $F$.  Because $S_k$ is isotropic for $\beta$, $G$ is
isotropic for $\beta_F$.  
By the universal property of the
isotropic flag variety, there exists a unique section
$\sigma:M\rightarrow C'$ such that the pullback of the universal flag
$S_2(F,\beta_F) \subset \pi^* F$ equals $G\subset F$.  On $C'$
there is a rank $k$ locally direct summand $S'_k$ of $\pi^* S_{k-2}^\perp$
defined as the preimage of the rank $2$ locally direct summand
$S_2(E,\beta_F)$ of  $\pi^*F = \pi^*(S_{k-2}^\perp/S_{k-2})$.
Because $S_2(F,\beta_F)$ is isotropic for $\beta_F$,
$S'_k$ is isotropic for $\beta$.  By the universal property of
$\text{Flag}_k(W,\beta)$, there exists a unique morphism
$g:C'\rightarrow \text{Flag}_k(W,\beta)$ such that $g^*S_k$ equals
$S'_k$ as a locally direct summand of $W\otimes_\kappa \OO_{C'}$.

\begin{prop} \label{prop-flagD5}
\marpar{prop-flagD5}
Assume $n$ equals $2k$. 
The morphism $\pi':C'\rightarrow M$ is a proper, smooth morphism, and
every fiber is a disjoint union of 2 copies of $\PP^1$.  There is a
unique open and closed subscheme $C$ containing the image of $\sigma$
such that $\pi:C\rightarrow M$ is a $\PP^1$-bundle.  The datum
$(\pi:C\rightarrow M,\sigma,g)$ is a family of pointed lines in
$\text{Flag}_k(W,\beta)$ parametrized by
$M$.  The associated morphism $\iota:M\rightarrow
\M_{0,1}(\text{Flag}_k(W,\beta),1)$ is an isomorphism, and
$\text{ev}\circ\iota$ is the tautological projection
$\text{Flag}_{k-2,k}(W,\beta) \rightarrow \text{Flag}_k(W,\beta)$.
\end{prop}

\begin{proof}
The first part follows from the fact that for $n=4$,
$\text{Flag}_2(W,\beta)$ is a disjoint union of two copies of
$\PP^1$.  Each $\PP^1$ has degree $2$ with respect to the Pl\"ucker
embedding, thus it has degree $1$ with respect to $\OO_X(1)$.
Therefore $(\pi,\sigma,g)$ is a family of pointed lines in $X$.  As in
the proofs of Proposition~\ref{prop-iota} and Proposition
~\ref{prop-flagBD}, to prove $\iota$ is an isomorphism, it suffices to
prove it is bijective on points.

\medskip\noindent
Let $C$ be a smooth conic in $\text{Flag}(k,W)$.
By ~\cite[Lemma 1]{Buch}, there is a rank $k-2$ subspace $S_{k-2}$ of $W$
such that every point of $C$
parametrizes a subspace $S_k$ containing $S_{k-2}$.
If $C$ is contained in $\text{Flag}_k(W,\beta)$, then $S_k$
is isotropic, hence also $S_{k-2}$ is isotropic.  Because $S_k$ is
isotropic and contains $S_{k-2}$, $S_k$ is contained in the rank $k+2$
annihilator, $S_{k-2}^\perp$.  Denote by $F$ the quotient space
$S_{k-2}^\perp/S_{k-2}$.  Associated to $\beta$ there is a symmetric
bilinear pairing $\beta_F$ on $F$.  A rank $k$ subspace $S_k$ of $W$
containing $S_{k-2}$ is isotropic for $\beta$ iff the subspace
$S_k/S_{k-2}$ of $F$ is isotropic for $\beta_F$.  Therefore $C$ is
also a smooth conic in $\text{Flag}_2(F,\beta_F)$.  As above,
$\text{Flag}_2(F,\beta_F)$ is a disjoint union of two smooth conics,
i.e., $C$ is one of the two connected components of
$\text{Flag}_2(F,\beta_F)$.  Also, since the subspaces of $F$
parametrized by $\text{Flag}_2(F,\beta_F)$ collectively span $F$ and
have common intersection $(0)$, $S_{k-2}$ is the common intersection in
$W$ of the spaces $S_k$ for every point $[S_k]$ of $C$.  Therefore the
space of conics in $\text{Flag}_k(W,\beta)$ is the bijective image 
of the \'etale double
cover $\text{Flag}_2(F,\beta_F)$ of $\text{Flag}_{k-2}(W,\beta)$.
Therefore $\iota$ is bijective.
\end{proof}

\begin{cor} \label{cor-flagD5}
\marpar{cor-flagD5/prop-flagD6}
Assume $n$ equals $2k$.  The pullbacks under $\iota$ of $T_\text{ev}$
and $\psi^\vee$ admit canonical isomorphisms,
$$
\iota^*T_{\text{ev}}\cong (S_k/S_{k-2})\otimes S_{k-2}^\vee,
$$
$$
\iota^*\psi^\vee \cong \bigwedge^2(S_k/S_{k-2})^\vee.
$$
\end{cor}

\begin{proof}
The projection $\text{Flag}_{k-2,k}(W,\beta) \rightarrow
\text{Flag}_k(W,\beta)$ is the relative Grassmannian
$\text{Flag}(k-2,S_k)$.  So the first isomorphism follows from the
well-known isomorphism of the vertical tangent bundle of a relative
Grassmannian.  By construction, $\pi:C \rightarrow M$ is one connected
component of the relative isotropic Grassmannian
$\text{Flag}_2(F,\beta_F)$, where $F := S_{k-2}^\perp/S_{k-2}$.  The
vertical tangent bundle of $\text{Flag}(2,F)$ equals
$S(2,F)^\vee\otimes \pi^*F/S(2,F)$.  The restriction of this sheaf to
$\text{Flag}_2(F,\beta_F)$ is $S_2(F,\beta_F)^\vee \otimes
S_2(F,\beta_F)^\vee$.  
By ~\cite[p. 474]{HarrisTu}, the
normal bundle of the regular embedding
$\text{Flag}_2(F,\beta_F)\rightarrow \text{Flag}(2,F)$ equals
$\text{Sym}^2(S_2(F,\beta_F))^\vee$.  Therefore the vertical tangent
bundle of $\text{Flag}_2(F,\beta_F)$ is the kernel of the map
$S_2(F,\beta_F)^\vee \otimes S_2(F,\beta_F)^\vee \rightarrow
\text{Sym}^2(S_2(F,\beta_F))^\vee$, i.e., $\bigwedge^2
S_2(F,\beta_F)^\vee$.  The bundle $\psi^\vee$ is the pullback of the
vertical tangent bundle by the section $\sigma$.  Since $\sigma$ pulls
back $S_2(F,\beta_F)$ to $G=S_k/S_{k-2}$, the bundle $\psi^\vee$ is
canonically isomorphic to
$\bigwedge^2(S_k/S_{k-2})^\vee$. 
\end{proof}

\medskip\noindent
\textbf{Skew-symmetric case, $n\geq 2k$.}
Assume $(W,\beta)$ is a skew-symmetric pairing of dimension $n$.
Assume $n$ is at least $2k$.  
There is a natural embedding of the isotropic Grassmannian in the
classical Grassmannian, $\text{Flag}_k(W,\beta)\subset
\text{Flag}(k,W)$.  Denote by $\OO_X(1)$ the pullback of
$\OO_{\PP(\bigwedge^k W)}(1)$ under the Pl\"ucker embedding.  

\medskip\noindent
Denote by $\M$ the space $\M_{0,1}(X,1)$.  
Denote by $M_\text{pre}$ the isotropic flag variety,
$$
M = \text{Flag}_{k-1,k}(W,\beta).
$$
Denote by $S_{k-1,\text{pre}}\subset S_{k,\text{pre}} \subset
W\otimes_\kappa \OO_{M_\text{pre}}$ the universal isotropic flag.
On $M_\text{pre}$ there is a rank $n-k+1$ locally direct summand
$S_{k-1,\text{pre}}^\perp$ of $W\otimes_\kappa \OO_{M_\text{pre}}$
defined as the annihilator of $S_{k-1,\text{pre}}$ under $\beta$.
Since $S_{k,\text{pre}}$ and $S_{k-1,\text{pre}}$ are isotropic, 
each one
is a locally direct summand of
$S_{k-1,\text{pre}}$.  Denote by $G$ the rank $n+1-2k$, locally free sheaf
$S_{k-1,\text{pre}}^\perp/S_{k,\text{pre}}$.  

\medskip\noindent
Denote by $\rho:M\rightarrow M_{\text{pre}}$ the projective bundle,
$$
M = \text{Flag}(1,G) = \PP(G).
$$
Denote by $S_{k-1}$, resp. $S_k$, the locally direct summand of
$W\otimes_\kappa \OO_M$ obtained by pulling back $S_{k-1,\text{pre}}$,
resp. $S_{k,\text{pre}}$.  Each one is an isotropic locally direct
summand of $W\otimes_\kappa \OO_M$.  Also the
annihilator
$S_{k-1}^\perp$ of $S_{k-1}$ with respect to $\beta$ equals
$\rho^*S_{k-1,\text{pre}}^\perp$ as subsheaves of $W\otimes_\kappa
\OO_M$.  
There is a rank $k+1$ locally direct
summand $R_{k+1}$ of $S_{k-1}^\perp$ defined as the preimage of the
rank $1$ locally direct summand $\OO_{PP(G)}(-1)$ of $\pi^*G =
S_{k-1}^\perp/S_k$.  

\medskip\noindent
Altogether, this defines a flag of locally direct summands of
$W\otimes_\kappa \OO_M$,
$$
S_{k-1}\subset S_k \subset R_{k+1} \subset W\otimes_\kappa \OO_M.
$$
The first two terms are isotropic.  The term $R_{k+1}$ is not
necessarily isotropic, but it is contained in $S_{k-1}^\perp$.  The
following lemma is straightforward.

\begin{lem} \label{lem-SSR}
\marpar{lem-SSR}
The flag $S_{k-1}\subset S_k\subset R_{k+1} \subset W\otimes_\kappa
\OO_M$ is the universal $(k-1,k,k+1)$-flag of locally direct summands
of $W$ such that $S_k$ is isotropic and $R_{k+1}$ is contained in
$S_{k-1}^\perp$.
\end{lem}

\begin{rmk} \label{rmk-SSR}
There is a natural action of $\textbf{Sp}(W,\beta)$ on $M$.  There are
2 orbits.  The closed orbit is the projective homogeneous space
$\text{Flag}_{k-1,k,k+1}(W,\beta)$.  The complement of the closed
orbit is a non-projective homogeneous space.  If $1<k<n/2$, the
automorphism group of $M$ equals $\textbf{Sp}(W,\beta)$, thus $M$ is
not a homogeneous space for any group.
\end{rmk}

\medskip\noindent
Define $E'$ to be the rank $2$ locally free sheaf, $R_{k+1}/S_{k-1}$,
and define $L'$ to be the invertible sheaf, $S_k/S_{k-1}$.  Denote by
$\pi:\PP(E') \rightarrow M$ the associated $\PP^1$-bundle.  There is a
unique section $\sigma':M\rightarrow \PP(E')$ such that the pullback of
the rank $1$ locally direct summand $\OO_{\PP(E')}(-1)$ of $\pi^* E'$
equals $L'$.

\medskip\noindent
On $\PP(E')$ there is a rank $k$ locally direct summand $S'_k$ 
of $\pi^* R_{k+1}$
defined as the preimage of the rank $1$
locally direct summand $\OO_{\PP(E')}(-1)$ of
$\pi^* E' = \pi^*(R_{k+1}/S_{k-1})$.  Because $\pi^*R_{k+1}$ is a rank
$k+1$ 
locally direct summand of $W\otimes_\kappa \OO_{\PP(E')}$, $S'_k$ is a
rank $k$
locally direct summand of $W\otimes_\kappa \OO_{\PP(E)}$.
Denote by $F$ the quotient $S_{k-1}^\perp/S_{k-1}$.    
Associated to $\beta$
there is a skew-symmetric pairing $\beta_F$ for $F$.  
The locally direct summand $\OO_{\PP(E')}(-1)$ of
$\pi^* F$ is isotropic for $\beta_F$ because
\emph{every} rank $1$ locally direct summand of a skew-symmetric
pairing is isotropic.  Because $S'_k/\pi^*S_{k-1})$ is isotropic for
$\beta_F$, $S'_k$ is isotropic for $\beta$.
By the
universal property of $\text{Flag}_k(W,\beta)$, there exists a unique
morphism $g':\PP(E')\rightarrow \text{Flag}_k(W,\beta)$ such that $(g')^*
S_k$ equals $S'_k$ as a locally direct summand of $W\otimes_\kappa
\OO_{\PP(E')}$.  

\begin{prop} \label{prop-flagC7}
\marpar{prop-flagC7}
Assume $n$ is at least $2k$.
The datum $(\pi',\sigma',g')$ is a family of pointed
lines in $\text{Flag}_k(W,\beta)$ parametrized by  $M$.  The associated
morphism $\iota:M\rightarrow \M$ is an
isomorphism. 
\end{prop}

\begin{proof}
As in the proof of Proposition~\ref{prop-flagBD},  the first statement
follows from Proposition~\ref{prop-iota}, and the second statement
reduces to surjectivity of $\iota$.
Let
$[S_{k-1}\subset S_k\subset R_k\subset W]$ be a flag parametrizing a
line in $\text{Flag}(k,W)$ contained in $\text{Flag}_k(W,\beta)$.
Then $S_k$ is isotropic, and hence also $S_{k-1}$ is isotropic.  Every
vector in $R_k$ is contained in a subspace $S'_k$ containing
$S_{k-1}$.  Because $S'_k$ is isotropic, $S'_k$ is contained in
$S_{k-1}^\perp$.  Therefore $R_k$ is contained in $S_{k-1}^\perp$.  By
Lemma~\ref{lem-SSR}, $[S_{k-1}\subset S_k\subset R_k\subset W]$ is
contained in the image of $\iota$.
\end{proof}

\begin{cor} \label{cor-flagC7}
\marpar{cor-flagC7/prop-flagC8}
Assume $n$ is at least $2k$.
The pullback under $\iota$ of $T_\text{ev}$ admits
a short exact sequence,
$$
0 \rightarrow [(R_{k+1}/S_k)^\vee\otimes
(S_{k-1}^\perp/R_{k+1})] \rightarrow \iota^* T_\text{ev} \rightarrow
[(S_k/S_{k-1})\otimes S_{k-1}^\vee] \rightarrow 0.  
$$ 
And the pullback of $\psi^\vee$ under $\iota$ admits an isomorphism,
$$
\iota^*\psi^\vee \cong (R_{k+1}/S_k)\otimes (S_k/S_{k-1})^\vee.
$$
\end{cor}

\begin{proof}
This is very similar to the proof of Corollary~\ref{cor-flagBD}.
\end{proof}

\section{Maps of vector bundles on the projective line}
\label{subsec-P1}
\marpar{subsec-P1}
Using the propositions of Subsections ~\ref{subsec-pl} and
~\ref{subsec-BCD1}, a very twisting family of pointed lines on a
classical or isotropic Grassmannian is equivalent to a flag of locally
direct summands of $W\otimes_\kappa \OO_{\PP^2}$ such that the
associated locally free sheaf $\zeta^*T_\text{ev}$ is ample and the
associated invertible sheaf
$\zeta^*\psi^\vee$ has nonnegative degree.  
However, in the isotropic case, the flags are difficult to construct
directly.  
This subsection contains the proof of a fact
about maps of vector bundles on the projective line which will be
useful for constructing flags.

\medskip\noindent
Let $a \leq b$ be nonnegative integers.  Define $H$ to be the rank
$ab(b+1-a)$ vector space,
$$
H =
\text{Hom}_{\OO_{\PP^1}}(\OO_{\PP^1}^{\oplus
  a},\OO_{\PP^1}(b-a)^{\oplus b}).
$$
There is an open subset of $H$ parametrizing maps whose cokernel is
locally free of rank $b-a$.  There is an open subset of this subset
parametrizing maps whose cokernel is isomorphic to
$\OO_{\PP^1}(b)^{\oplus(b-a)}$.  Denote by $H^o$ this open subset of
$H$.

\begin{prop} \label{prop-H}
\marpar{prop-H}
The open subset $H^o$ is not empty.
\end{prop}

\begin{proof}
If $a$ is zero, this is vacuous.  Thus assume $a\geq 1$.
Let $U$ be a rank $2$ vector space and identify $\PP^1$ with
$\PP(U)$.  The homogeneous coordinate ring of $\PP^1$ is
$S:= \text{Sym}^\bullet(U^\vee)$.  Denote the graded pieces by $S_k$,
i.e., $S_k = \text{Sym}^k(U^\vee)$.  By convention, define $S_{-1}$ to
be the zero vector space. 
Denote by $A$ the associative, unital $\kappa$-algebra of linear maps
from $S$ to $S$, $A=\text{Hom}(S,S)$.  This has a natural structure of
left $S$-module by $(p\cdot L)(q) = pL(q)$ for every $p,q$ in $S$ and
every $L$ in $A$.
Denote by 
$\text{Diff}$ the subalgebra of $\text{Hom}(S,S)$ of
differential operators on $S$.  This is a left $S$-submodule of $A$.

\medskip\noindent  
Denote by $d$ the unique $\kappa$-derivation,
$$
d:S \rightarrow
S \otimes_\kappa U^\vee,
$$
such that $d|_{S_1}: U^\vee \rightarrow S\otimes_\kappa U^\vee$
factors as
$$
U^\vee \xrightarrow{\text{Id}} \kappa\otimes_\kappa U^\vee =
S_0\otimes_\kappa U^\vee \subset S\otimes_\kappa U^\vee.
$$
The derivation $d$ identifies $U$ as a
linear subspace of $\text{Diff}$.
Define
$\text{Diff}_{a-1}$ to be the linear subspace of $\text{Diff}$ generated
by $U^{\otimes (a-1)}$.  This, of course, is just an isomorphic copy of
$\text{Sym}^{a-1}(U)$.  In particular it has rank $a$.  Also every
element in $\text{Diff}_{a-1}$ has order $a-1$.  

\medskip\noindent
There is a canonical map $c:A \rightarrow
(\bigwedge^2 U)\otimes_\kappa A$ defined as follows.  Choose an ordered
basis $\mathbf{e}_0, \mathbf{e}_1$ for $U$, and denote by $T_0,T_1$
the dual ordered basis for $U^\vee$.  For every linear map
$L:S\rightarrow S$, define $c(L):S\rightarrow (\bigwedge^2 U)\otimes
S$
to be,
$$
c(L)(p) = (\mathbf{e}_0\wedge\mathbf{e}_1)\otimes
(T_1 L(T_0\cdot p) - T_0 L(T_1\cdot p)).
$$
It is straightforward to check
this is independent of the choice of basis.  The map $c$ is a morphism
of left $S$-modules (where $\bigwedge^2 U$ is given the trivial
$S$-module structure).    
More importantly, if $L$ is a differential operator of order $k+1$,
then $c(L)$ is a differential operator of order $k$.  Define
$c^l:A \rightarrow (\bigwedge^2 U)^{\otimes l}\otimes_\kappa A$ to be
the $l$-fold composition of $c$ in the obvious way.
Then $\text{Ker}(c^l)$ contains the subspace of differential operators of
order $\leq l-1$.  

\medskip\noindent
For every integer $k$,
define $Q_k$ to be the quotient $\text{Hom}(S_k,S)$ of $A$, and denote
by $\pi_k:A\rightarrow Q_k$ the quotient map.  By convention, define
$Q_{-1}$ to be the zero vector space.  
The space $Q_k$ has 
a natural left $S$-module structure by $(q\cdot L)(p) = qL(p)$ and
$\pi_k$ is a morphism of left $S$-modules.  
Make $Q_k$
a graded $S$-module by defining $(Q_k)_l =
\text{Hom}(S_k,S_{k+l})$ for every integer $l$.  
The associated sheaf on $\PP(U)$ is
$S_k^\vee\otimes_\kappa \OO_{\PP(U)}(k)$.  
For every integer $0\leq l \leq k+1$, there is a
unique degree $2l$ map of graded $S$-modules, $C_l:Q_k \rightarrow
(\bigwedge^2 U)^{\otimes l}\otimes_\kappa Q_{k-l}$, such that the
following diagram commutes,
$$
\begin{CD}
A @> c^l >> (\bigwedge^2 U)^{\otimes l} \otimes_\kappa A \\
@V \pi_k VV @VV \text{Id}\otimes \pi_{k-l} V \\
Q_k @> C_l >> (\bigwedge^2 U)^{\otimes l} \otimes_\kappa Q_{k-l}.
\end{CD}
$$
This induces a map of associated sheaves, 
$$
\widetilde{C}_l:
S_k^\vee \otimes_\kappa
\OO_{\PP(U)}(k) \rightarrow [(\bigwedge^2 U)^{\otimes l}\otimes_\kappa
  S_{k-l}^\vee]\otimes_\kappa \OO_{\PP(U)}(k+l).
$$

\medskip\noindent
The composite map $\text{Diff}_{a-1} \hookrightarrow A
\xrightarrow{\pi_{b-1}} Q_{b-1}$ 
has image in the subspace $(Q_{b-1})_{-(a-1)}$, and so
induces a map of associated sheaves,
$$
\phi_{a,b}:\text{Diff}_{a-1}\otimes_\kappa \OO_{\PP(U)} \rightarrow
S_{b-1}^\vee \otimes_\kappa \OO_{\PP(U)}(b-a).
$$
Twisting the map $\widetilde{C}_a$ appropriately gives a map of
associated sheaves,
$$
\psi_{a,b}:S_{b-1}^\vee \otimes_\kappa [(\bigwedge^2 U)^{\otimes l} \otimes
  S_{b-a-1}^\vee] \otimes_\kappa \OO_{\PP(U)}(b).
$$
Since $\text{dim}(\text{Diff}_{a-1})$ equals $a$,
$\text{dim}(S_{b-1}^\vee)$ equals $b$, and $\text{dim}((\bigwedge^2
U)^{\otimes l}\otimes S_{b-a-1}^\vee)$ equals $b-a$, the proposition
is implied by the following.

\begin{claim} \label{claim-ses}
\marpar{claim-ses}
The following sequence of sheaves on $\OO_{\PP(U)}$ is exact,
$$
0 \rightarrow \text{Diff}_{a-1}\otimes_\kappa \OO_{\PP(U)}
\xrightarrow{\phi_{a,b}} S_{b-1}^\vee \otimes_\kappa \OO_{\PP(U)}(b-a)
 \xrightarrow{\psi_{a,b}} [(\bigwedge^2 U)^\vee \otimes
   S_{a+b-1}^\vee] \otimes_\kappa \OO_{\PP(U)}(b) \rightarrow 0.
$$
\end{claim}

\medskip\noindent
First of all, because $\text{Diff}_{a-1}$ is contained in the kernel
of $c^a$, 
$\psi_{a,b}\circ \phi_{a,b}$ is the zero map.  There is a natural
action of $\textbf{GL}(U)$ on $\PP(U)$.  
Each
vector bundle in the claim
has a natural $\mathbf{GL}(U)$-linearization and each of 
$\phi_{a,b}$ and $\psi_{a,b}$ is $\mathbf{GL}(U)$-equivariant.  Thus
to prove the 
sequence is exact, it suffices to prove it is exact at one point of
$\PP(U)$.  With respect to the bases $\mathbf{e}_0,\mathbf{e}_1$ of
$U$ and the dual basis $T_0,T_1$, denote by $\partial_0$,
$\partial_1$ the elements of $\text{Diff}_1$ corresponding to
$\mathbf{e}_0,\mathbf{e}_1$.  Then an ordered basis for
$\text{Diff}_{a-1}$ consists of,
$$
\frac{1}{(a-1)!\cdot 0!}\partial_0^{a-1}, \dots,
\frac{1}{(a-k)!\cdot (k-1)!}\partial_0^{a-k}\partial_1^{k-1}, \dots,
\frac{1}{0!\cdot (a-1)!} \partial_1^{a-1}.
$$
An ordered basis for $S_{b-1}$ consists of,
$$
T_0^{b-1}, \dots, T_0^{b-j}T_1^{j-1}, \cdots, T_1^{b-1},
$$
and this gives a dual ordered basis for $S_{b-1}^\vee$.  Similarly, an
ordered basis for $S_{b-a-1}$ consists of,
$$
T_0^{a+b-1}, \dots, T_0^{a+b-1-i}T_1^{i-1},\dots, T_1^{a+b-1},
$$
and, tensored with $(\textbf{e}_0\wedge\textbf{e}_1)^{\otimes l}$, 
this gives a dual ordered basis for $(\bigwedge^2 U)^{\otimes
  l}\otimes S_{a+b-1}^\vee$.  With
respect to these ordered bases, the entries of the
matrix of $\phi_{a,b}$ equal,
$$
(\phi_{a,b})_{j,k} = \binom{b-j}{a-k} \binom{j-1}{k-1} T_0^{b-1+k-j}
T_1^{j-k},
$$
and the entries of the matrix of $\psi_{a,b}$ equal,
$$
(\psi_{a,b})_{i,j} = (-1)^{j-i} \binom{a}{j-i} T_0^{j-i} T_1^{a+i-j}.
$$
Plugging in $T_0=1, T_1=0$ gives the matrices,
$$
\phi_{a,b}|_{[1,0]} = \lt(
\begin{array}{c}
D_{a,a} \\
\hline
0_{b-a,a}
\end{array} \rt),
$$
and
$$
\psi_{a,b}|_{[1,0]} = (-1)^a\cdot \lt(
\begin{array}{c|c}
0_{a,b-a} & I_{b-a,b-a}
\end{array} \rt),
$$ 
where $0_{k,l}$ is the zero $k\times l$ matrix, $I_{k,k}$ is the
$k\times k$ identity matrix, and 
$D_{a,a}$ is the invertible $a\times a$ diagonal matrix with entries
$(D_{a,a})_{j,j} =
\binom{b-j}{a-j}$.  
These matrices visibly give a short exact
sequence of $\kappa$-vector spaces.  
\end{proof}

\medskip\noindent
Because the entries of $\phi_{a,b}$ and $\psi_{a,b}$ are monomials,
the very twisting families constructed later are
more likely to be orbit curves.

\section{Application to isotropic subspaces} \label{subsec-app}
\marpar{subsec-app}
This subsection applies Proposition~\ref{prop-H} to construct some
useful isotropic subspaces of $W\otimes_\kappa \OO_{\PP^1}$.  

\begin{hyp} \label{hyp-app}
\marpar{hyp-app}
Let
$a$ and $b$ be positive integers.  Assume that $b$ is at least $2a$.  
\end{hyp}

\medskip\noindent
Let $W_{b,+}$
be an $b$-dimensional $\kappa$-vector space.  
Denote the dual vector space by $W_{b,-}$.  
By Proposition
~\ref{prop-H}, there exists a map,
$$
\phi_{a,b,+}:\OO_{\PP^1}(-(b-a))^{\oplus a} \rightarrow
W_{b,+}\otimes_\kappa \OO_{\PP^1},
$$
whose cokernel is isomorphic to $\OO_{\PP^1}(a)^{\oplus (b-a)}$.
Thus the annihilator in $W_{b,-} \otimes_\kappa \OO_{\PP^1}$
of the image of $\phi_{a,b,+}$ is the image of a map,
$$
\psi_{a,b,+}^\dagger:\OO_{\PP^1}(-a)^{\oplus (b-a)} \rightarrow
W_{b,-}\otimes_\kappa \OO_{\PP^1}.  
$$
By hypothesis, $b-a$ is at least $a$.
Thus, by Proposition
~\ref{prop-H}, there exists a map,
$$
\phi_{a,b-a}:\OO_{\PP^1}(-(b-a))^{\oplus a} \rightarrow
\OO_{\PP^1}(-a)^{\oplus (b-a)},
$$
whose cokernel is isomorphic to $\OO_{\PP^1}^{\oplus(b-2a)}$.
Define $\phi_{a,b,-}$ to be the composite map,
$$
\psi_{a,b,+}^\dagger \circ \phi_{a,b-a}:\OO_{\PP^1}(-(b-a))^{\oplus
  a} \rightarrow W_{b,-} \otimes_\kappa \OO_{\PP^1}.
$$

\medskip\noindent
Define $W_{2b}$ to be $W_{b,+}\oplus W_{b,-}$.  Thus
$(W_{2b})^\vee$ is canonically isomorphic to $W_{b,-}\oplus
W_{b,+}$.  
In the symmetric case, define
$\beta_{2b}: W_{2b} \rightarrow W_{2b}^\vee$ to be the linear map
$(W_{b,+}\oplus W_{b,-}) \rightarrow (W_{b,-} \oplus
W_{b,+})$ with matrix, 
$$
\beta_{2b} = \lt(
\begin{array}{c|c}
0 & \text{Id}_{W_{b,-}} \\
\hline
\text{Id}_{W_{b,+}} & 0 
\end{array} \rt).
$$
In the skew-symmetric case, define $\beta_{2b}$ to be the linear map with
matrix, 
$$
\beta_{2b} = \lt(
\begin{array}{c|c}
0 & -\text{Id}_{W_{b,-}} \\
\hline
\text{Id}_{W_{b,+}} & 0 
\end{array} \rt).
$$   
Define $E_{2a,2b}$ to be the image in $W_{2b}\otimes_\kappa
\OO_{\PP^1}$
of the sheaf map,
$$
\phi_{a,b,+}\oplus\phi_{a,b,-}:\OO_{\PP^1}(-(b-a))^{\oplus a}
\oplus \OO_{\PP^1}(-(b-a))^{\oplus a} \rightarrow
(W_{b,+}\otimes_\kappa \OO_{\PP^1}) \oplus
(W_{b,-}\otimes_\kappa \OO_{\PP^1}). 
$$
Denote by $E_{2a,2b}^\perp$ the annihilator of $E_{2a,2b}$ for
$\beta_{2b}$.  

\begin{lem} \label{lem-app}
\marpar{lem-app}
The subsheaf $E_{2a,2b}$ is a
rank $2a$ locally direct summand of $W_{2b}\otimes_\kappa
\OO_{\PP^1}$ isotropic for $\beta_{2b}$.  
The quotient
$(E_{2a,2b})^\perp/E_{2a,2b}$ is isomorphic to
$\OO_{\PP^1}^{\oplus (2b-4a)}$.  
And, $E_{2a,2b}^\vee$ is an ample
vector bundle on $\PP^1$.
\end{lem}

\begin{proof}
Denote by $F$ the rank $2b-4a)$ locally free sheaf
$E_{2a,2b}^\perp/E_{2a,2b}$.  Associated to $\beta$, there
is a symmetric, resp. skew-symmetric, pairing $\beta_F$ for $F$.  By
construction,
$G: =\text{Image}(\psi_{a,b,+}^\dagger)/\text{Image}(\phi_{a,b,-})$ is an
rank $b-2a$ locally direct summand of $F$, isotropic for $\beta_F$.
Therefore $F$ is isomorphic to $G\oplus G^\vee$.  By construction, $G$
is isomorphic to $\OO_{\PP^1}^{\oplus (b-2a)}$.  Therefore $F$ is
isomorphic to $\OO_{\PP^1}^{\oplus (2b-4a)}$.  

\medskip\noindent
By construction, $E_{2a,2b}^\vee$ is isomorphic to
$\OO_{\PP^1}(b-a)^{\oplus 2a}$.  By hypothesis, $a$ is positive, and
$b-2a$ is nonnegative, so also $b-a = (b-2a)+a \geq a >0$.  Thus
$\OO_{\PP^1}(b-a)$ is ample.  Since $2a > 0$,
$\OO_{\PP^1}(b-a)^{\oplus 2a}$ is ample.
\end{proof}

\section{The classical Grassmannian} \label{subsec-A}
\marpar{subsec-A}
Let $n>2$ be an integer and let $k$ be an integer $0<k<n$.  Let
$V$ be an $n$-dimensional $\kappa$-vector space, and denote by
$(X,\OO_X(1))$
the Grassmannian $\text{Flag}(k,V)$ and the Pl\"ucker invertible
sheaf.  
Replacing $\text{Flag}(k,V)$ by the isomorphic scheme
$\text{Flag}(n-k,V^\vee)$ if necessary, assume $k\leq n/2$.  Of
course $\text{Flag}(k,V)$ equals $\mathbf{SL}_n/P_k$, where $P_k$ is the
maximal parabolic group corresponding to the $k^\text{th}$ node of the
Dynkin diagram $A_{n-1}$.

\medskip\noindent
By Proposition ~\ref{prop-iota},
a morphism $\zeta:\PP^1 \rightarrow
\M_{0,1}(\text{Flag}(k,V),1)$ is equivalent to a $(k-1,k,k+1)$-flag of 
locally direct summands,
$$
E_{k-1}\subset E_k \subset E_{k+1}\subset V\otimes_\kappa \OO_{\PP^1}.
$$
By Corollary ~\ref{cor-iota}, the morphism $\zeta$ is very
twisting iff 
\begin{enumerate}
\item[(i)]
the bundle,
$$
[(E_{k+1}/E_k)^\vee \otimes ((V\otimes \OO_{\PP^1})/E_{k+1})]\oplus
[(E_k/E_{k-1})\otimes E_{k-1}^\vee],
$$
is ample, and 
\item[(ii)]
the bundle,
$$
(E_{k+1}/E_k)\otimes (E_k/E_{k-1})^\vee,
$$
has nonnegative degree.
\end{enumerate}

\begin{prop} \label{prop-grasslines2}
\marpar{prop-grasslines2}
There exists a very twisting morphism $\zeta:\PP^1\rightarrow
\M_{0,1}(X,1)$.
\end{prop}

\begin{proof}
It is equivalent to prove there
exists a flag $E_{k-1}\subset E_k \subset E_{k+1}\subset V\otimes
\OO_{\PP^1}$ satisfying Conditions (i)--(ii).  Let $T_0,T_1$
denote homogeneous coordinates on $\PP^1$.

\medskip\noindent
Let $V'_{2k-2}$ be a rank $2k-2$ vector space, let $E'_{k-1}$ be
$\OO_{\PP^1}(-1)^{\oplus(k-1)}$, and let
$\phi'_{k-1,2k-2}:E'_{k-1}\rightarrow V'_{2k-2}\otimes_\kappa
\OO_{\PP^1}$ be a morphism whose cokernel is isomorphic to
$\OO_{\PP^1}(1)^{\oplus (k-1)}$, as in Proposition ~\ref{prop-H}.  Let
$V''_{n+1-2k}$ be a rank $n+1-2k$ vector space, let $E''_1$ be
$\OO_{\PP^1}(-(n-2k))$, and let $\phi''_{1,n+1-2k}:E''_1\rightarrow
V''_{n+1-2k}\otimes_\kappa \OO_{\PP^1}$ be a morphism whose cokernel
is isomorphic to $\OO_{\PP^1}(1)^{n-2k}$, as in
Proposition~\ref{prop-H}.  Finally, let $V'''_1$ be a rank $1$ vector
space, let $E'''_1$ be $V'''_1\otimes_\kappa \OO_{\PP^1}$, and let
$\phi'''_{1,1}$ be the identity map.  

\medskip\noindent
Define $V$ to be the direct sum of $V'_{2k-2}$, $V''_{n+1-2k}$ and
$V'''_1$.  Define $E_{k+1}$ to be the image of $\phi'_{k-1,2k-1}\oplus
\phi''_{1,n+1-2k}\oplus \phi'''_{1,1}$ in $V\otimes \OO_{\PP^1}$.  
The cokernel is
isomorphic to $\OO_{\PP^1}(1)^{\oplus(k-1)}\oplus
\OO_{\PP^1}(1)^{\oplus(n-2k)}$, i.e.,
$\OO_{\PP^1}(1)^{\oplus(n-k-1)}$. 

\medskip\noindent
Define $E_k$ to be the image of $\phi'_{k-1,2k-2}\oplus
\phi''_{1,n+1-2k}$.  
The quotient
$E_{k+1}/E_k$ equals $E'''_{1} \cong \OO_{\PP^1}$.  In particular,
$(E_{k+1}/E_k)^\vee\otimes ((V\otimes \OO_{\PP^1})/E_{k+1})$ is
isomorphic to $\OO_{\PP^1}(1)^{n-k-1}$, which is ample.  This is half of
Condition (i).
For $E_{k-1}$,
there are two cases.

\medskip\noindent
\textbf{Case I: $k=1$.}
In this case, define $E_{k-1}$ to be $(0)$.  The quotient
$E_k/E_{k-1}$ equals $E''_{1} \cong
\OO_{\PP^1}(-(n-2k))$.  So $\text{deg}(E_{k+1}/E_k)$ equals $0$ and
$\text{deg}(E_k/E_{k-1})$ equals $-(n-2k)$.  Since $n\geq 2k$,
$-(n-2k)\leq 0$, i.e., Condition (ii) holds.  Since $k=1$, Condition
(i) holds.

\medskip\noindent
\textbf{Case II: $k > 1$.}
Decompose $E'_{k-1}$ as $E'_{k-2,a}\oplus E'_{1,b}$, where
$E'_{k-2,a}$ is the first $k-2$ summands and $E'_{1,b}$ is the last
summand.  
Define $E''''_{1}$ to be $\OO_{\PP^1}(-(n+1-2k))$.  Define
$\phi'_{k-2,a}:E'_{k-2,a} \rightarrow E'_{k-2,a}$ to be the identity
map.  And define $\phi''''_{1,2}:E''''_1\rightarrow
E'_{1,b} \oplus E''_1$, i.e., $\OO_{\PP^1}(-(n+1-2k)) \rightarrow
\OO_{\PP^1}(-1)\oplus \OO_{\PP^1}(-(n-2k))$, to be the map with
matrix,
$$
\phi''''_{1,2} = \lt(
\begin{array}{l}
T_0^{n-2k} \\
T_1
\end{array}
\rt)
$$
Define $E_{k-1}$ to be $E'_{k-2,a}\oplus E''''_1$, and define
$\phi_{k-1}:E_{k-1}\rightarrow E_k$ to be
$\phi'_{k-2,a}\oplus \phi''''_{1,2}$.
The map $\phi''''_{1,2}$ is injective with cokernel $\OO_{\PP^1}$.  Thus
$\phi_{k-1}$ is injective and
$E_k/E_{k-1}$ equals $\OO_{\PP^1}$.  Since both
$E_{k+1}/E_k$ and $E_k/E_{k-1}$ have degree $0$, Condition (ii) holds.
Also, $E_{k-1}^\vee \otimes (E_k/E_{k-1})$ equals
$\OO_{\PP^1}(1)^{\oplus(k-2)}\oplus \OO_{\PP^1}(n+1-2k)$.  Since
$n\geq 2k$, this is an ample bundle.  Therefore Condition (i) holds.
\end{proof}

\begin{claim} \label{claim-grasslines}
\marpar{claim-grasslines}
The very twisting family in the proof of
Proposition~\ref{prop-grasslines2} can be chosen to be an orbit curve.
\end{claim}

\begin{proof}
For simplicity, assume $k>1$; the case $k=1$ is similar and easier.
Choose $\phi'_{k-1,2k-2}$ to be the map with matrix, 
$$
\phi'_{k-1,2k-2} = \lt(
\begin{array}{cccc}
T_0 & 0   & \dots & 0 \\
T_1 & 0   & \dots & 0 \\
0   & T_0 & \dots & 0 \\
0   & T_1 & \dots & 0 \\
\vdots & \vdots & \ddots & \vdots \\
0   & 0   & \dots & T_0 \\
0   & 0   & \dots & T_1 
\end{array} \rt)
$$
Choose $\phi''_{1,n+1-2k}$ to be the map with matrix,
$$
\phi''_{1,n+1-2k} = \lt(
\begin{array}{c}
T_0^{n-2k} \\
\vdots \\
T_0^{n-2k-j}T_1^j \\
\vdots \\
T_1^{n-2k}
\end{array}
\rt)
$$

\medskip\noindent
Let $\lambda':\mathbf{G}_m\times
V\rightarrow V$ be the linear action compatible with the direct sum
decomposition $V'_{2k-2}\oplus V''_{n+1-2k}\oplus V'''_1$ 
given by the diagonal matrices
$D_1$, $D_2$,
$$
D_1 = \lt(
\begin{array}{cc|cc|c|cc}
1 & 0 & 0 & 0 & \dots & 0 & 0 \\
0 & t & 0 & 0 & \dots & 0 & 0 \\
\hline
0 & 0 & 1 & 0 & \dots & 0 & 0 \\
0 & 0 & 0 & t & \dots & 0 & 0 \\
\hline
\vdots & \vdots & \vdots & \vdots & \ddots & \vdots & \vdots \\
\hline
0 & 0 & 0 & 0 & \dots & 1 & 0 \\
0 & 0 & 0 & 0 & \dots & 0 & t
\end{array}
\rt)
$$
$$
D_2 = \lt(
\begin{array}{cccc}
t &   0 & \dots  & 0 \\
0 & t^2 & \dots  & 0 \\
\vdots & \vdots  & \ddots & \vdots \\
0 &   0 &  \dots & t^{n+1-2k}
\end{array}
\rt)
$$
and $D_3$ is the $1\times 1$ matrix $t^c$ for
$c=-[(k-1)+(n+2-2k)(n+1-2k)/2]$.  Let $E'_{k-1}$ be the subspace of
$V'_{2k-2}$ which is the image of the matrix,
$$
\lt( 
\begin{array}{cccc}
1 & 0 & \dots & 0 \\
1 & 0 & \dots & 0 \\
0 & 1 & \dots & 0 \\
0 & 1 & \dots & 0 \\
\vdots & \vdots & \ddots & \vdots \\
0 & 0 & \dots & 1 \\
0 & 0 & \dots & 1 
\end{array}
\rt).
$$
Let $E''_{1}$ be the subspace of $V''_{n+1-2k}$ spanned by the vector,
$$
\lt(
\begin{array}{c}
1 \\
\vdots \\
1
\end{array}
\rt).
$$
And let $E'''_{1}$ equals $V'''_1$.  Define $E_{k+1} = E'_{k-1}\oplus
E''_{1}\oplus E'''_{1}$.  Define $E_k$ to be $E'_{k-1}\oplus
E''_{1}$.

\medskip\noindent
Decompose $V'_{2k-2}$ as $V'_{2k-4,a}\oplus V'_{2,b}$ where 
$V'_{2k-4,a}$ is the first $2k-4$
summands, and $V'_{2,b}$ is the last 2 summands.
Define $E'_{k-2,a}$ to be the subspace of $V'_{2k-4,a}$ which is the image of
the matrix,
$$
\lt( 
\begin{array}{cccc}
1 & 0 & \dots & 0 \\
1 & 0 & \dots & 0 \\
0 & 1 & \dots & 0 \\
0 & 1 & \dots & 0 \\
\vdots & \vdots & \ddots & \vdots \\
0 & 0 & \dots & 1 \\
0 & 0 & \dots & 1 
\end{array}
\rt).
$$
Define $E''''_{1}$ to be the subspace of $V'_{2,b}\oplus V''_{n+1-2k}$
spanned by 
the vector,
$$
\lt( 
\begin{array}{c}
1 \\
\vdots \\
1 \\
\hline
1 \\
\vdots \\
1
\end{array}
\rt).
$$
Define $E_{k-1}$ to be $E'_{k-2,a}\oplus E''''_{1}$.  This gives a flag
of subbundle of $V$,
$E_{k-1}\subset E_k\subset E_{k+1}\subset V$. 
Define $P$ to be the maximal parabolic subgroup of
$\mathbf{SL}(V)$ that is the
stabilizer of the flag $E_k\subset V$
and define $P'$ to be the stabilizer of the flag $E_{k-1}\subset E_k
\subset E_{k+1}\subset V$.  The 1-parameter subgroup 
$\lambda:\mathbf{G}_m\rightarrow
\mathbf{SL}(V)$ is defined above.  
The rational
curve $\zeta:\PP^1 \rightarrow \text{Flag}(k-1,k,k+1,V)$ equals the
orbit curve associated to $\lambda$ and the flag $[E_{k-1}\subset
  E_k\subset E_{k+1}\subset V]$.
\end{proof}

\section{Isotropic Grassmannians, Case I}
\label{subsec-keven}
\marpar{subsec-keven}
Let $(W,\beta)$ be a symmetric or skew-symmetric pairing of dimension
$n$.  If $n=2$ or $3$
there is no very twisting family of pointed lines to an
isotropic Grassmannian of $(W,\beta)$.
Thus assume $n\geq 4$.  
This subsection proves existence of a very twisting family of pointed
lines on the Grassmannian of isotropic $k$-planes
when $k$ is odd and $n\geq 2k+2$.  

\begin{hyp} \label{hyp-caseI}
\marpar{hyp-caseI}
The pairing $(W,\beta)$ is symmetric of dimension $2m$ or $2m+1$ or
the pairing $(W,\beta)$ is skew-symmetric of dimension $2m$, and
$k=2l+1$.  Assume that $l$ is nonnegative and $m$ is at least $\max(2,2l+2)$.
\end{hyp}

\medskip\noindent
The last inequality is equivalent to $k$ is positive, $n\geq 4$, and
$n\geq 2k+2$. 

\medskip\noindent
Let $(W'_4,\beta_4)$ and $E'_{2,4}$ be as in Subsection ~\ref{subsec-app} for
$a=1$ and $b=2$.  In particular, $E'_{2,4}$ is isomorphic to
$\OO_{\PP^1}(-1)^{\oplus 2}$.  

\medskip\noindent
\textbf{Case Ia: $k=1$.}  
Assume that $k=1$.
Let $(W''_{n-4},\beta'')$ be a
symmetric pairing, resp. skew-symmetric pairing, 
of dimension $n-4$.  Define $(W,\beta)$ to be the orthogonal direct sum
of $(W'_4,\beta'_4)$ and $(W''_{n-4},\beta'')$.  Define $E_2$ to be
$E'_{2,4}$ 
considered as a locally direct summand of $W\otimes_\kappa
\OO_{\PP^1}$ via the embedding $W'_4\otimes_\kappa
\OO_{\PP^1}\hookrightarrow W\otimes_\kappa \OO_{\PP^1}$.  
Define $E_1$ to be a direct summand
$\OO_{\PP^1}(-1)$ in $E_2$, and define $E_0$ to be the zero sheaf.

\begin{lem} \label{lem-caseIa}
\marpar{lem-caseIa}
Assume $k=1$ and $n\geq 4$.
The flag $E_0\subset E_1\subset E_2\subset W\otimes_\kappa
\OO_{\PP^1}$ is a $(0,1,2)$-flag of isotropic locally direct summands
for $\beta$.  
The cokernel $E_2^\perp/E_2$ is
isomorphic to $W''_{n-4}\otimes_\kappa
\OO_{\PP^1}$. 
The cokernels $E_2/E_1$ and
$E_1/E_0$ are isomorphic to $\OO_{\PP^1}(-1)$.  And $E_0^\vee$ is the
zero sheaf.
\end{lem}

\begin{proof}
By construction, $E_2$ is isotropic of rank $2$.  Since $E_1$ is
contained in $E_2$, it is isotropic.  Of course $(0)$ is isotropic.
By construction, the annihilator $E_2^\perp$ of $E_2$ with respect to
$\beta$
equals the direct sum of the annihilator $(E'_{2,4})^\perp$ of
$E'_{2,4}$ with respect to $\beta'_{2,4}$ and
$W''_{n-4}\otimes_\kappa\OO_{\PP^1}$.  By Lemma~\ref{lem-app},
$(E'_{2,4})^\perp$ equals $E'_{2,4}$.  Therefore $E_2^\perp/E_2$
equals $W''_{n-4}\otimes_\kappa \OO_{\PP^1}$. 
By definition of $E_1$ and $E_0$, $E_2/E_1 \cong E_1/E_0 \cong
\OO_{\PP^1}(-1)$.  The dual of the zero sheaf is the zero sheaf. 
\end{proof}

\begin{prop} \label{prop-caseIa}
\marpar{prop-caseIa}
Assume $k=1$.  In the skew-symmetric case, assume $n\geq 4$.  In the
symmetric case, assume $n\geq 5$.  
The morphism $\zeta:\PP^1 \rightarrow M$ associated to the flag in
Lemma~\ref{lem-caseIa} is a very twisting family of pointed lines on
$\text{Flag}_1(W,\beta)$.  
\end{prop}

\begin{proof}
By Corollary~\ref{cor-flagBD} and Corollary ~\ref{cor-flagC7}, in both
the symmetric and skew-symmetric case $\iota^*\psi^\vee$ equals
$(E_{k+1}/E_k)\otimes (E_k/E_{k-1})^\vee$.  By Lemma~\ref{lem-caseIa},
this equals $\OO_{\PP^1}$, and so has nonnegative degree.  This is
(iii) of Definition ~\ref{defn-vt}.  In the skew-symmetric case, by
Corollary~\ref{cor-flagC7} and since $E_0^\vee$ is the zero sheaf, 
$\iota^* T_\text{ev}$ is isomorphic to
$(E_{2}/E_1)^\vee \otimes (E_0^\perp/E_2)$. 
Also since $E_0$ is
the zero sheaf, $E_0^\perp$ equals $W\otimes_\kappa \OO_{\PP^1}$.  In
particular, $E_0^\perp/E_2$ has rank $n-2$, which is positive by the
hypothesis that $n \geq 4$.  
Since $W\otimes_\kappa \OO_{\PP^1}$ is globally generated, the
quotient $E_0^\perp/E_2$ is globally generated of positive rank.  The
tensor product of an ample bundle and a globally generated, positive
rank bundle is an ample bundle.   
Since
$(E_2/E_1)^\vee$ is ample, the tensor product $(E_2/E_1)^\vee \otimes
(E_0^\perp/E_2)$ is ample.  

\medskip\noindent
The argument in the symmetric case is the same, except $E_0^\perp/E_2$
is replaced by $E_2^\perp/E_2$, which is globally generated of rank
$n-4$ by Lemma~\ref{lem-caseIa}.  
The rank $n-4$ is positive by the hypothesis that $n\geq 5$.
\end{proof}

\medskip\noindent
An argument similar to the proof of Claim~\ref{claim-grasslines}
proves the very twisting family can be chosen to be an orbit curve.

\medskip\noindent
\textbf{Case Ib: $k>1$.}
Assume now that $k>1$, i.e., $l\geq 1$.  
Let $(W''_{2m-4},\beta''_{2m-4})$ and $E''_{2l,2m-4,\text{pre}}$ be as in
Subsection ~\ref{subsec-app} for $a=l$ and $b=m-2$.
Hypothesis~\ref{hyp-caseI} 
implies $b\geq 2a$, i.e., Hypothesis~\ref{hyp-app} holds.  
By Lemma~\ref{lem-app}, $(E''_{2l,2m-4,\text{pre}})^\vee$ is ample.
Let $f:\PP^1 \rightarrow \PP^1$ be any finite morphism such that
$f^*[(E''_{2l,2m-4,\text{pre}})^\vee] \otimes \OO_{\PP^1}(-1)$ is
ample.  In every case except $(l,m)=(1,4)$, it suffices to take $f$ to
be the identity map.  If $(l,m)=(1,4)$, it suffices to take $f$ to be
any finite morphism of degree $\geq 2$.  At any rate, define
$E''_{2l,2m-4}$ to be $f^*(E''_{2l,2m-4,\text{pre}})$ considered
as a subsheaf of $f^*(W''_{2m-4}\otimes_\kappa \OO_{\PP^1}) =
W''_{2m-4}\otimes_\kappa \OO_{\PP^1}$.  

\medskip\noindent
If $n=2m$, define $(W,\beta)$ to be the orthogonal direct sum of
$(W'_4,\beta'_4)$ and $(W''_{2m-4},\beta''_{2m-4})$.  If $n=2m+1$,
which can only occur in the symmetric case,
let $(\mathbf{1},\beta_{\mathbf{1}})$ be a
symmetric pairing of dimension $1$, and 
define $(W,\beta)$ to
be the orthogonal direct sum of $(W'_4,\beta'_4)$,
$(W''_{2m-4},\beta''_{2m-4})$ and $(\mathbf{1},\beta_{\mathbf{1}})$.
Define $E_{2l+2}$ to be the direct sum $E'_{2,4}$ and 
$E''_{2l,2m-4}$.  
Define $E_{2l+1}$ to be the direct sum of one
direct summand $\OO_{\PP^1}(-1)$ of $E'_{2,4}$ and 
$E''_{2l,2m-4}$.  Finally, define $E_{2l}$ to be $E''_{2l,2m-4}$.

\begin{lem} \label{lem-caseIb'}
\marpar{lem-caseIb'}
Assume $l\geq 1$ and $m\geq 2l+2$.  The flag $E_{2l}\subset E_{2l+1}
\subset E_{2l+2} \subset W\otimes_\kappa \OO_{\PP^1}$ is a
$(k-1,k,k+1)$-flag of isotropic locally direct summands for $\beta$.
The cokernel $E_{2l+2}^\perp/E_{2l+2}$ is isomorphic to
$\OO_{\PP^1}^{\oplus (2m-4l-4)}$ if $n=2m$, 
respectively $\OO_{\PP^1}^{\oplus
  (2m-4l-3)}$ if 
$n=2m+1$.  The cokernel $E_{2l}^\perp/E_{2l+2}$ is isomorphic to
$\OO_{\PP^1}(1)^{\oplus 2} \oplus \OO_{\PP^1}^{\oplus (2m-4l-4)}$ if
$n=2m$, 
respectively $\OO_{\PP^1}(1)^{\oplus 2} \oplus \OO_{\PP^1}^{\oplus
  (2m-4l-3)}$ if $n=2m+1$.  The cokernels $E_{2l+2}/E_{2l+1}$ and
$E_{2l+1}/E_{2l}$ are each isomorphic to $\OO_{\PP^1}(-1)$.  And
$E_{2l}^\vee \otimes (E_{2l+1}/E_{2l})$ is ample.
\end{lem}

\begin{proof}
Since $E'_{2,4}$ is isotropic for $\beta'_{2,4}$
and $E''_{2l,2m-4}$ is isotropic for $\beta''_{2l,2m-4}$, $E_{2l+2}$
is isotropic for $\beta$.  Since $E_{2l+1}$ and $E_{2l}$ are contained
in $E_{2l+2}$, they are also isotropic for $\beta$.  The annihilator
$E_{2l+2}^\perp$ of $E_{2l+2}$ with respect to $\beta$ is the direct
sum of the annihilator $(E'_{2,4})^\perp$ of $E'_{2,4}$ with respect
to $\beta'_{2,4}$, the annihilator $(E''_{2l,2m-4})^\perp$ of
$E''_{2l,2m-4}$ with respect to $\beta''_{2l,2m-4}$, and also
$\OO_{\PP^1}$ if $n=2m+1$.  Therefore $E_{2l+2}^\perp/E_{2l+2}$ equals
the direct sum of $(E'_{2,4})^\perp/E'_{2,4}$,
$(E''_{2l,2m-4})^\perp/E''_{2l,2m-4}$, and also $\OO_{\PP^1}$ if
$n=2m+1$.   By Lemma~\ref{lem-app}, $E'_{2,4}$ is its own annihilator,
and $(E''_{2l,2m-4})^\perp/E''_{2l,2m-4}$ equals $\OO_{\PP^1}^{\oplus
  (2m-4l-4)}$.  
Therefore $E_{2l+2}^\perp/E_{2l+2}$ equals
$\OO_{\PP^1}^{\oplus (2m-4l-4)}$ if $n=2m$, and equals
$\OO_{\PP^1}^{\oplus (2m-4l-3)}$ if $n=2m+1$.  

\medskip\noindent
The computation of $E_{2l}^\perp/E_{2l+2}$ is similar, except the
summand $(E'_{2,4})^\perp/E'_{2,4}$ is replaced by $W'_4\otimes_\kappa
\OO_{\PP^1}/E'_{2,4}$.  Since $E_{2,4}'$ equals its own annihilator,
$W'_4\otimes_\kappa \OO_{\PP^1}/E'_{2,4}$ equals the dual of
$E'_{2,4}$.  Thus $W'_4\otimes_\kappa \OO_{\PP^1}/E'_{2,4}$ equals
$(\OO_{\PP^1}(-1)^{\oplus 2})^\vee$, i.e., $\OO_{\PP^1}(1)^{\oplus
  2}$.  Therefore $E_{2l}^\perp/E_{2l+2}$ equals
$\OO_{\PP^1}(1)^{\oplus 2} \oplus \OO_{\PP^1}^{\oplus (2m-4l-4)}$ if
$n=2m$, and $\OO_{\PP^1}(1)^{\oplus 2}\oplus \OO_{\PP^1}^{\oplus
  (2m-4l-3)}$ if $n=2m+1$.

\medskip\noindent
By definition, $E_{2l+2}/E_{2l+1}$ and $E_{2l+1}/E_{2l}$ are
isomorphic to $\OO_{\PP^1}(-1)$ and
$E_{2l}^\vee\otimes(E_{2l+1}/E_{2l})$, i.e., $E_{2l}^\vee \otimes
\OO_{\PP^1}(-1)$, is ample.
\end{proof}

\begin{prop} \label{prop-caseIb}
\marpar{prop-caseIb}
Assume $l\geq 1$ and $m\geq 2l+1$. 
The morphism $\zeta:\PP^1 \rightarrow M$ associated to the flag in
Lemma~\ref{lem-caseIb'} is a very twisting family of pointed lines on
$\text{Flag}_k(W,\beta)$. 
\end{prop}

\begin{proof}
This is very similar to the proof of Proposition ~\ref{prop-caseIa}.
\end{proof}

\medskip\noindent
It seems likely $\zeta$ can be chosen to be an orbit curve.  However,
since the entries of the matrix for $\phi''_{2l,2m-4}$ are typically
not monomials, it is not certain.  

\section{Isotropic Grassmannians, Case II}
\label{subsec-caseII}
\marpar{subsec-caseII}
Let $(W,\beta)$ be a symmetric or skew-symmetric pairing of dimension $n$.
This subsection proves existence of a very twisting family of pointed
lines on the Grassmannian of isotropic $k$-planes when $k$ is even and $n\geq
2k+2$.  

\begin{hyp} \label{hyp-caseII}
\marpar{hyp-caseII}
The pairing $(W,\beta)$ is symmetric of dimension $2m$ or $2m+1$ or
the pairing $(W,\beta)$ is skew-symmetric of dimension $2m$, and
$k=2l$.  Assume that $l$ is positive and $m$ is at least $2l+1$.
\end{hyp}

\medskip\noindent
The last inequality is equivalent to $k$ is positive and $n\geq
2k+2$.  In particular, observe that $m$ is at least $3$.

\medskip\noindent
Define $(W'_6,\beta'_6)$ as in Subsection ~\ref{subsec-app} for
$b=3$. 
Choose an
ordered basis $\textbf{e}_1,\textbf{e}_2,\textbf{e}_3$ of $W'_{3,+}$.
Denote the dual ordered basis of $W'_{3,-}$ by $x_1, x_2, x_3$.  
Define
$E'_{3,6}$ to be the rank $3$ locally free $\OO_{\PP^1}$-module
$\OO_{\PP^1}(-2)\textbf{f} \oplus \OO_{\PP^1}(-1)\textbf{g}_1 \oplus
\OO_{\PP^1}(-1)\textbf{g}_2$, where the symbols $\textbf{f},
\textbf{g}_1,\textbf{g}_2$ are simply place-holders.  There is an
isomorphism $\phi'_{3,6}$ of $E'_{3,6}$ to an isotropic locally direct
summand of $W'_6\otimes_\kappa \OO_{\PP^1}$.  The definition of
$\phi'_{3,6}$ is different in the symmetric and skew-symmetric case.

\medskip\noindent
\textbf{Symmetric case.}
Define $\phi'_{3,6}:E'_{3,6} \rightarrow W'_6\otimes_\kappa
\OO_{\PP^1}$ to be the unique $\OO_{\PP^1}$-module homomorphism
satisfying, 
$$
\begin{array}{ccc}
\phi'_{3,6}(\textbf{g})   & = & T_0^2\textbf{e}_1 + T_0T_1\textbf{e}_2 +
T_1^2 x_3 \\
\phi'_{3,6}(\textbf{f}_1) & = & T_0 \textbf{e}_3 - T_1 x_2 \\
\phi'_{3,6}(\textbf{f}_2) & = & T_1 x_1 - T_0 x_2
\end{array}
$$
The proof of the following lemma is a straightforward computation.

\begin{lem} \label{lem-caseIIa'}
\marpar{lem-caseIIa'}
Let $(W'_6,\beta'_6)$ be the symmetric pairing from above.  
The image of $\phi'_{3,6}$ is a rank $3$ locally direct summand of
$W'_6\otimes_\kappa \OO_{\PP^1}$ isotropic for $\beta'_6$.  It equals
its own annihilator with respect to $\beta'_6$.  
\end{lem}

\medskip\noindent
\textbf{Skew-symmetric case.}
Define $\phi'_{3,5}:E'_{3,6} \rightarrow W'_6\otimes_\kappa
\OO_{\PP^1}$ to be the unique $\OO_{\PP^1}$-module homomorphism
satisfying, 
$$
\begin{array}{ccc}
\phi'_{3,6}(\textbf{g})   & = & -T_0^2\textbf{e}_1 +
T_0T_1(\textbf{e}_2-x_2) + T_1^2 x_3 \\
\phi'_{3,6}(\textbf{f}_1) & = & T_0 (\textbf{e}_2+x_2) + 2T_1 x_1 \\
\phi'_{3,6}(\textbf{f}_2) & = & T_1 (\textbf{e}_2+x_2) + 2T_0 \textbf{e}_3
\end{array}
$$
The proof of the following lemma is a straightforward computation.

\begin{lem} \label{lem-caseIIb'}
\marpar{lem-caseIIb'}
Let $(W'_6,\beta'_6)$ be the skew-symmetric pairing from above.
The image of $\phi'_{3,6}$ is a rank $3$ locally direct summand of
$W'_6\otimes_\kappa \OO_{\PP^1}$ isotropic for $\beta'_6$.  It equals
its own annihilator with respect to $\beta'_6$.  
\end{lem}

\medskip\noindent
\textbf{Case IIa: $k=2$.}
Assume that $k=2$, i.e., $l=1$.  Then Hypothesis ~\ref{hyp-caseII} is
equivalent to $n\geq 6$.  
Let $(W''_{n-6},\beta'')$ be a
symmetric pairing, resp. skew-symmetric pairing, 
of dimension $n-6$.  Define $(W,\beta)$ to be the orthogonal direct sum
of $(W'_6,\beta'_6)$ and $(W''_{n-6},\beta'')$.  Define $E_3$ to be
$E'_{3,6}$. 
Define $E_2$ to be direct summand
$\OO_{\PP^1}(-2)\textbf{g}\oplus \OO_{\PP^1}(-1)\textbf{f}_1$, 
and define $E_1$ to be $\OO_{\PP^1}(-2)\textbf{g}$.

\begin{lem} \label{lem-caseIIa}
\marpar{lem-caseIIa}
Assume $k=2$ and $n\geq 6$, either the symmetric or the skew-symmetric
case.  
The flag $E_1\subset E_2\subset E_3\subset W\otimes_\kappa
\OO_{\PP^1}$ is a $(1,2,3)$-flag of isotropic locally direct summands
for $\beta$.  
The cokernel $E_3^\perp/E_3$ is
isomorphic to $W'_{n-6}\otimes_\kappa \OO_{\PP^1}$.  The cokernel
$E_1^\perp/E_3$ is isomorphic to $\OO_{\PP^1}(1)^{\oplus 2}\oplus
W'_{n-4}\otimes_\kappa \OO_{\PP^1}$.  The cokernels $E_3/E_2$ and
$E_2/E_1$ are isomorphic to $\OO_{\PP^1}(-1)$.  And $E_1^\vee\otimes
(E_2/E_1)$ is isomorphic to the ample invertible sheaf
$\OO_{\PP^1}(1)$.  
\end{lem}

\begin{proof}
By Lemma~\ref{lem-caseIIa'}, resp. Lemma~\ref{lem-caseIIb'}, $E_3$ is
isotropic of rank $3$.  
Since $E_2$ and $E_1$ are contained in $E_3$, they are
isotropic.  By construction, 
the annihilator $E_3^\perp$ of $E_3$ with respect to
$\beta$ is the direct sum of the annihilator $(E'_{3,6})^\perp$ of
$E'_{3,6}$ with respect to $\beta'_{3,6}$ and $W_{n-6}''\otimes_\kappa
\OO_{\PP^1}$.  By Lemma~\ref{lem-caseIIa'},
resp. Lemma~\ref{lem-caseIIb'}, $(E'_{3,6})^\perp$
equals $E'_{3,6}$.  Therefore $E_3^\perp/E_3$ equals
$W_{n-6}''\otimes_\kappa \OO_{\PP^1}$.  Similarly, $E_1^\perp$ is the
direct sum of the annihilator of $\OO_{\PP^1}(-2)\textbf{g}$ with
respect to $\beta'_{3,6}$ and $W_{n-6}''\otimes_\kappa \OO_{\PP^1}$.
Since $E'_{3,6}$ equals its own annihilator,
$(\OO_{\PP^1}(-2)\textbf{g})^\perp/E'_{3,6}$ equals the dual of
$E'_{3,6}/\OO_{\PP^1}(-2)\textbf{g}$.  Thus
$(\OO_{\PP^1}(-2)\textbf{g})^\perp/E'_{3,6}$ equals $(\OO_{\PP^1}(-1)^{\oplus
  2})^\vee$, i.e., $\OO_{\PP^1}(1)^{\oplus 2}$.  Therefore $E_1^\perp/E_3$
equals $\OO_{\PP^1}(1)^{\oplus 2}\oplus W_{n-6}''\otimes_\kappa
\OO_{\PP^1}$.  

\medskip\noindent
By the definition of $E_2$ and $E_1$, $E_3/E_2 \cong E_2/E_1 \cong
\OO_{\PP^1}(-1)$.  Since $E_1\cong \OO_{\PP^1}(-2)$, $E_1^\vee
\otimes(E_2/E_1)$ equals $\OO_{\PP^1}(2)\otimes \OO_{\PP^1}(-1)$,
i.e., $\OO_{\PP^1}(1)$.  
\end{proof}

\begin{prop} \label{prop-caseIIa}
\marpar{prop-caseIa}
Assume $k=2$ and $n\geq 6$, either the symmetric or the skew-symmetric
case.  
The morphism $\zeta:\PP^1 \rightarrow M$ associated to the flag in
Lemma~\ref{lem-caseIIa} is a very twisting family of pointed lines on
$\text{Flag}_2(W,\beta)$.  
\end{prop}

\begin{proof}
This is very similar to the proof of Proposition ~\ref{prop-caseIa}.
Of course now the term $E_1^\vee \otimes (E_2/E_2)$ in $\zeta^*
T_\text{ev}$ is nonzero.  But by the last part of
Lemma~\ref{lem-caseIIa}, this is ample. 
\end{proof}

\medskip\noindent
An argument similar to the proof of Claim~\ref{claim-grasslines}
proves the very twisting family can be chosen to be an orbit curve.

\medskip\noindent
\textbf{Case IIb: $k>2$.}
Assume now that $k>2$, i.e., $l\geq 2$.
In both the symmetric case and the skew-symmetric case,
let $(W''_{2m-6},\beta''_{2m-6})$ and $E''_{2l-2,2m-6,\text{pre}}$ be as in
Subsection ~\ref{subsec-app} for $a=l-1$ and $b=m-3$.  Because $l\geq
2$, $a$ is positive.  By Hypothesis~\ref{hyp-caseII}, $b\geq 2a$.  In
particular, $b$ is positive and  Hypothesis~\ref{hyp-app} holds.
By Lemma~\ref{lem-app}, $(E''_{2l-2,2m-6},\text{pre})^\vee$ is ample.
Let $f:\PP^1 \rightarrow \PP^1$ be any finite morphism such that
$f^*[(E''_{2l-2,2m-6,\text{pre}})^\vee] \otimes \OO_{\PP^1}(-1)$ is
ample.  In every case except $(l,m)=(2,5)$, it suffices to take $f$ to
be the identity map.  If $(l,m)=(2,5)$, it suffices to take $f$ to be
any finite morphism of degree $\geq 2$.  At any rate, define
$E''_{2l-2,2m-6}$ to be $f^*(E''_{2l-2,2m-6,\text{pre}})$ considered
as a subsheaf of $f^*(W''_{2m-6}\otimes_\kappa \OO_{\PP^1}) =
W''_{2m-6}\otimes_\kappa \OO_{\PP^1}$.  

\medskip\noindent
If $n=2m$, define $(W,\beta)$ to be the orthogonal direct sum of
$(W'_6,\beta'_6)$ and $(W''_{2m-6},\beta''_{2m-6})$.  If $n=2m+1$,
which can only occur in the symmetric case,
let $(\mathbf{1},\beta_{\mathbf{1}})$ be a
symmetric pairing of dimension $1$, and 
define $(W,\beta)$ to
be the orthogonal direct sum of $(W'_6,\beta'_6)$,
$(W''_{2m-6},\beta''_{2m-6})$ and $(\mathbf{1},\beta_{\mathbf{1}})$.
Define $E_{2l+1}$ to be the direct sum $E'_{3,6}$ and 
$E''_{2l-2,2m-6}$.  
Define $E_{2l}$ to be the direct sum of
$\OO_{\PP^1}(-2)\textbf{g}\oplus \OO_{\PP^1}(-1)\textbf{f}_1$ and
$E''_{2l-2,2m-6}$.  Finally, define $E_{2l-1}$ to be the direct sum of
$\OO_{\PP^1}(-2)\textbf{g}$ and 
$E''_{2l-2,2m-6}$.

\begin{lem} \label{lem-caseIIb}
\marpar{lem-caseIIb}
Assume $l\geq 2$ and $m\geq 2l+1$.  The flag $E_{2l-1}\subset E_{2l}
\subset E_{2l+1} \subset W\otimes_\kappa \OO_{\PP^1}$ is a
$(k-1,k,k+1)$-flag of isotropic locally direct summands for $\beta$.
The cokernel $E_{2l+1}^\perp/E_{2l+1}$ is isomorphic to
$\OO_{\PP^1}^{\oplus (2m-4l-2)}$ if $n=2m$, 
respectively $\OO_{\PP^1}^{\oplus
  (2m-4l-1)}$ if 
$n=2m+1$.  The cokernel $E_{2l-1}^\perp/E_{2l+1}$ is isomorphic to
$\OO_{\PP^1}(1)^{\oplus 2} \oplus \OO_{\PP^1}^{\oplus (2m-4l-2)}$ if
$n=2m$, 
respectively $\OO_{\PP^1}(1)^{\oplus 2} \oplus \OO_{\PP^1}^{\oplus
  (2m-4l-1)}$ if $n=2m+1$.  The cokernels $E_{2l+1}/E_{2l}$ and
$E_{2l}/E_{2l-1}$ are each isomorphic to $\OO_{\PP^1}(-1)$.  And
$E_{2l-1}^\vee \otimes (E_{2l}/E_{2l-1})$ is ample.
\end{lem}

\begin{proof}
Since $E_{3,6}'$ is isotropic for $\beta_{3,6}'$ and $E_{2l-2,2m-6}''$
is isotropic for $\beta_{2l-2,2m-6}''$, $E_{2l+1}$ is isotropic for
$\beta$.  Since $E_{2l}$ and $E{2l-1}$ are contained in $E_{2l+1}$,
they are also isotropic for $\beta$.  The annihilator $E_{2l+1}^\perp$
of $E_{2l+1}$ with respect to $\beta$ is the direct sum of the
annihilator $(E_{3,6}')^\perp$ of $E_{3,6}'$ with respect to
$\beta_{3,6}'$, the annihilator $(E_{2l-2,2m-6}'')^\perp$ of
$E_{2l-2,2m-6}''$ with respect to $\beta_{2l-2,2m-6}''$, and also
$\OO_{\PP^1}$ if $n=2m+1$.  Therefore $E_{2l+1}^\perp/E_{2l+2}$ equals
the direct sum of $(E'_{3,6})^\perp/E'_{3,6}$,
$(E''_{2l-2,2m-6})^\perp/E''_{2l-2,2m-6}$, and also $\OO_{\PP^1}$ if
$n=2m+1$.  By Lemma ~\ref{lem-caseIIa'}, resp. Lemma
~\ref{lem-caseIIb'}, $E'_{3,6}$ is its own annihilator.  By Lemma
~\ref{lem-app}, $(E''_{2l-2,2m-6})^\perp/E_{2l-2,2m-6}''$ equals
$\OO_{\PP^1}^{\oplus (2m-4l-2)}$.  
Therefore $E_{2l+1}^\perp/E_{2l+1}$ equals
$\OO_{\PP^1}^{\oplus (2m-4l-2)}$ if $n=2m$, and equals
$\OO_{\PP^1}^{\oplus (2m-4l-1)}$ if $n=2m+1$.

\medskip\noindent
The computation of $E_{2l-1}^\perp/E_{2l+1}$ is similar, except the
summand $(E'_{3,6})^\perp/E'_{3,6}$ is replaced by
$(\OO_{\PP^1}(-2)\textbf{g})^\perp /E'_{3,6}$.  
Since $E_{3,6}'$ equals its own annihilator,
$(\OO_{\PP^1}(-2)\textbf{g})^\perp /E'_{3,6}$
equals the dual of
$E'_{3,6}/\OO_{\PP^1}(-2)\textbf{g}$.  
Thus 
$(\OO_{\PP^1}(-2)\textbf{g})^\perp /E'_{3,6}$
equals
$(\OO_{\PP^1}(-1)^{\oplus 2})^\vee$, i.e., $\OO_{\PP^1}(1)^{\oplus
  2}$.  Therefore $E_{2l-1}^\perp/E_{2l+1}$ equals
$\OO_{\PP^1}(1)^{\oplus 2} \oplus \OO_{\PP^1}^{\oplus (2m-4l-2)}$ if
$n=2m$, and $\OO_{\PP^1}(1)^{\oplus 2}\oplus \OO_{\PP^1}^{\oplus
  (2m-4l-1)}$ if $n=2m+1$.

\medskip\noindent
By definition, $E_{2l+1}/E_{2l}$ and $E_{2l}/E_{2l-1}$ are
isomorphic to $\OO_{\PP^1}(-1)$.  Since $(E_{2l-2,2m-6}'')^\vee\otimes
\OO_{\PP^1}(-1)$ is ample, also
$E_{2l-1}^\vee\otimes(E_{2l}/E_{2l-1})$, which equals the direct sum of
$\OO_{\PP^1}(1)$ and $(E_{2l-2,2m-6}'')^\vee \otimes 
\OO_{\PP^1}(-1)$, is ample.
\end{proof}

\begin{prop} \label{prop-caseIIb}
\marpar{prop-caseIb}
Assume $l\geq 1$ and $m\geq 2l+1$. 
The morphism $\zeta:\PP^1 \rightarrow M$ associated to the flag in
Lemma~\ref{lem-caseIIb} is a very twisting family of pointed lines on
$\text{Flag}_k(W,\beta)$. 
\end{prop}

\begin{proof}
The proof is very similar to the proof of Proposition~\ref{prop-caseIIa}.
\end{proof}

\medskip\noindent
It seems likely $\zeta$ can be chosen to be an orbit curve.  However,
since the entries of the matrix for $\phi''_{2l,2m-4}$ are typically
not monomials, it is not certain.

\section{Isotropic Grassmannians, Case III}
\label{subsec-caseIII}
\marpar{subsec-caseIII}
Let $(W,\beta)$ be a symmetric pairing of dimension $n=2k$.  This
subsection proves existence of a very twisting family of pointed lines
on the Grassmannian of isotropic $k$-planes when $k\geq 3$.  There is
no very twisting family if $k=1$ or $k=2$.

\medskip\noindent
\textbf{Case IIIa, $k$ even.}
Let $l\geq 2$ be an integer, and let $k$ equal $2l$.  Let $n$ equal
$2k$, i.e., $4l$.  Let $(W'_4,\beta'_4)$ and $E'_{2,4}$ be as in
Subsection ~\ref{subsec-app} for $a=1$ and $b=2$.  In particular,
$E'_{2,4}$ is isomorphic to $\OO_{\PP^1}(-1)^{\oplus 2}$.  

\medskip\noindent
Let $(W''_{4l-4},\beta''_{4l-4})$ and $E''_{2l-2,4l-4,\text{pre}}$ be as in
Subsection ~\ref{subsec-app} for $a=l-1$ and $b=2l-2$.  Since $l\geq
2$, $a$ is positive.  And, of course, $b$ equals $2a$.
Thus
Hypothesis~\ref{hyp-app} holds.  
By Lemma~\ref{lem-app}, $(E''_{2l-2,4l-4,\text{pre}})^\vee$ is ample.
Let $f:\PP^1 \rightarrow \PP^1$ be any finite morphism such that
$f^*[(E''_{2l-2,4l-4,\text{pre}})^\vee] \otimes \OO_{\PP^1}(-1)$ is
ample.  In every case except $l=2$, it suffices to take $f$ to
be the identity map.  If $l=2$, it suffices to take $f$ to be
any finite morphism of degree $\geq 2$.  At any rate, define
$E''_{2l-2,4l-4}$ to be $f^*(E''_{2l-2,4l-4,\text{pre}})$ considered
as a subsheaf of $f^*(W''_{4l-4}\otimes_\kappa \OO_{\PP^1}) =
W''_{4l-4}\otimes_\kappa \OO_{\PP^1}$.  

\medskip\noindent
Define $(W,\beta)$ to be the orthogonal direct sum of
$(W'_4,\beta'_4)$ and $(W''_{4l-4},\beta''_{4l-4})$.  Define $E_{2l}$
to be the direct sum of $E'_{2,4}$ and $E''_{2l-2,4l-4}$.  And define
$E_{2l-2}$ to be $E''_{2l-2,4l-4}$.

\begin{lem} \label{lem-caseIIIa}
\marpar{lem-caseIIIb}
Assume $l\geq 2$.  Let $k$ equal $2l$ and let $n$ equal $4l$.  The
flag $E_{2l-2}\subset E_{2l} \subset W\otimes_\kappa \OO_{\PP^1}$ is a
$(k-2,k)$-flag of isotropic locally direct summands for $\beta$.  The
cokernel $E_{2l}/E_{2l-2}$ is isomorphic to $\OO_{\PP^1}(-1)^{\oplus
  2}$.  The determinant 
$\bigwedge^2(E_{2l}/E_{2l-2})^\vee$ is isomorphic to
$\OO_{\PP^1}(2)$.  And $(E_{2l}/E_{2l-2})\otimes E_{2l-2}^\vee$ is
ample.
\end{lem}

\begin{proof}
The proof is very similar to the proof of Lemma ~\ref{lem-caseIb'}.
The novel feature is the computation of
$\bigwedge^2(E_{2l}/E_{2l-2})^\vee$, which is obviously
$\OO_{\PP^1}(2)$ since $E_{2l}/E_{2l-2}$ equals
$\OO_{\PP^1}(-1)^{\oplus 2}$.  
\end{proof}

\begin{prop} \label{prop-caseIIIa}
\marpar{prop-caseIIIb}
Assume $l\geq 2$.  Let $k$ equal $2l$ and let $n$ equal $4l$.   
The morphism $\zeta:\PP^1 \rightarrow M$ associated to the flag in
Lemma~\ref{lem-caseIIIa} is a very twisting family of pointed lines on
$\text{Flag}_k(W,\beta)$.  
\end{prop}

\begin{proof}
The proof is very similar to the proof of
Proposition~\ref{prop-caseIa}.  Instead of
Corollary~\ref{cor-flagBD}, use Corollary~\ref{cor-flagD5}.
\end{proof}

\medskip\noindent
An argument similar to the proof of Claim~\ref{claim-grasslines}
proves the very twisting family can be chosen to be an orbit curve.

\medskip\noindent
\textbf{Case IIIb, $k$ odd.}
Let $l\geq 1$ be an integer, let $k$ equal $2l+1$, and let $n$ equal
$2k=4l+2$.  For $(n,k)=(2,1)$, there is no very twisting family of
lines.  Let $(W'_6,\beta'_6)$ and $E'_{3,6}$ be as in
Subsection~\ref{subsec-caseII} for the symmetric pairing.  In
particular, $E'_{3,6}$ is isomorphic to
$\OO_{\PP^1}(-2)\textbf{g}\oplus \OO_{\PP^1}(-1)\textbf{f}_1 \oplus
\OO_{\PP^1}(-1) \textbf{f}_2$.  If $l=1$, i.e., $k=3$ and $n=6$,
define $(W,\beta)$ to be $(W'_6,\beta'_6)$, define $E'_k$ to be
$E'_{3,6}$ and define $E'_{k-2}$ to be $\OO_{\PP^1}(-2)\textbf{g}$.  

\medskip\noindent
Next assume $l\geq 2$.
Let $(W''_{4l-4},\beta''_{4l-4})$ and $E''_{2l-2,4l-4,\text{pre}}$ be as in
Subsection ~\ref{subsec-app} for $a=l-1$ and $b=2l-2$.  Since $l\geq
2$, $a$ is positive.  And, of course, $b=2a$.
Thus
Hypothesis~\ref{hyp-app} holds.  
By Lemma~\ref{lem-app}, $(E''_{2l-2,4l-4,\text{pre}})^\vee$ is ample.
Let $f:\PP^1 \rightarrow \PP^1$ be any finite morphism such that
$f^*[(E''_{2l-2,4l-4,\text{pre}})^\vee] \otimes \OO_{\PP^1}(-1)$ is
ample.  In every case except $l=2$, it suffices to take $f$ to
be the identity map.  If $l=2$, it suffices to take $f$ to be
any finite morphism of degree $\geq 2$.  At any rate, define
$E''_{2l-2,4l-4}$ to be $f^*(E''_{2l-2,4l-4,\text{pre}})$ considered
as a subsheaf of $f^*(W''_{4l-4}\otimes_\kappa \OO_{\PP^1}) =
W''_{4l-4}\otimes_\kappa \OO_{\PP^1}$.

\medskip\noindent
Define $(W,\beta)$ to be the orthogonal direct sum of
$(W'_6,\beta'_6)$ and $(W''_{4l-4},\beta''_{4l-4})$.  Define $E_{2l+1}$
to be the direct sum of $E'_{3,6}$ and $E''_{2l-2,4l-4}$.  And define
$E_{2l-1}$ to be the direct sum of $\OO_{\PP^1}(-2)\textbf{g}$ and
$E''_{2l-2,4l-4}$.

\begin{lem} \label{lem-caseIIIb}
\marpar{lem-caseIIIb}
Assume $l\geq 1$.  Let $k$ equal $2l+1$ and let $n$ equal $4l+2$.  The
flag $E_{2l-1}\subset E_{2l+1} \subset W\otimes_\kappa \OO_{\PP^1}$ is a
$(k-2,k)$-flag of isotropic locally direct summands for $\beta$.  The
cokernel $E_{2l+1}/E_{2l-1}$ is isomorphic to $\OO_{\PP^1}(-1)^{\oplus
  2}$.  The determinant 
$\bigwedge^2(E_{2l+1}/E_{2l-1})^\vee$ is isomorphic to
$\OO_{\PP^1}(2)$.  And $(E_{2l+1}/E_{2l-1})\otimes E_{2l-1}^\vee$ is
ample.
\end{lem}

\begin{proof}
The proof is very similar to the proof of Lemma~\ref{lem-caseIIIa}.
\end{proof}

\begin{prop} \label{prop-caseIIIb}
\marpar{prop-caseIIIb}
Assume $l\geq 1$.  Let $k$ equal $2l+1$ and let $n$ equal $4l+2$.   
The morphism $\zeta:\PP^1 \rightarrow M$ associated to the flag in
Lemma~\ref{lem-caseIIIa} is a very twisting family of pointed lines on
$\text{Flag}_k(W,\beta)$.  
\end{prop}

\begin{proof}
The proof is very similar to the proof of Proposition ~\ref{prop-caseIIIa}.
\end{proof}

\medskip\noindent
An argument similar to the proof of Claim~\ref{claim-grasslines}
proves the very twisting family can be chosen to be an orbit curve.

\section{Isotropic Grassmannians, Case IV}
\label{subsec-caseIV}
\marpar{subsec-caseIV}
Let $(W,\beta)$ be a skew-symmetric pairing of dimension $n=2k$.  This
subsection proves existence of a very twisting family of pointed lines
on the Grassmannian of isotropic $k$-planes when $k\geq 2$.  For $k=1$ 
there is no very twisting family of pointed lines.  

\medskip\noindent
\textbf{Case IVa, $k$ even.}
Let $l\geq 1$ be an integer, and let $k$ equal $2l$.  Let $n$ equal
$2k$, i.e., $4l$.  Let $(W'_4,\beta'_2)$ be $(W'_{2,+}\oplus
W'_{2,-},\beta'_4)$ as in Subsection~\ref{subsec-app} for $b=2$.
Choose an ordered basis $\textbf{e}_1,\textbf{e}_2$ for $W'_{2,+}$,
and let $x_1,x_2$ be the dual ordered basis for $W'_{2,-}$.  Define
$R'_3$ to be the image of the sheaf,
$$
\phi'_3:\OO_{\PP^1}\{\textbf{a}\} \oplus
\OO_{\PP^1}(-1)\{\textbf{b}_+,\textbf{b}_-\} \rightarrow W'_4\otimes
\OO_{\PP^1},
$$
$$
\phi'_3(\textbf{a}) = \textbf{e}_2 - x_2,\ \phi'_3(\textbf{b}_+) =
T_0\textbf{e}_1 +T_1\textbf{e}_2, \ \phi'_3(\textbf{b}_-) = -T_1x_1 +
T_0x_2.
$$
Define $E'_2$ and $E'_1$ to be the subsheaves of $R'_3$ given by,
$$
\begin{array}{c}
E'_2 = \phi'_3(\OO_{\PP^1}(-1)\{\textbf{b}_+,\textbf{b}_-\}), \\
E'_1 = \phi'_3(\OO_{\PP^1}(-2)\{T_0\textbf{b}_+-T_1\textbf{b}_-\}) =
\OO_{\PP^1}(-2)\{T_0^2\textbf{e}_1 +T_0T_1(\textbf{e}_2-x_2) + T_1^2x_2\}.
\end{array}
$$

\begin{lem} \label{lem-caseIVa'}
\marpar{lem-caseIVa'}
Let $(W'_4,\beta'_4)$ be the skew-symmetric pairing from above.
The sheaf $E'_2$ is a rank $2$ isotropic locally direct summand of
$W'_4\otimes_\kappa \OO_{\PP^1}$.  The subsheaf $E'_1$ is a rank $1$
isotropic locally direct summand.  The annihilator of $E'_1$ with
respect to $\beta'_4$ equals $R'_3$.  The cokernels $R'_3/E'_2$ and
$E'_2/E'_1$ are both isomorphic to $\OO_{\PP^1}$.  And
$(E'_1)^\vee$ is an ample invertible sheaf.
\end{lem}
  
\begin{proof}
Since $\textbf{b}_+$ maps into $W'_{2,+}\otimes_\kappa \OO_{\PP^1}$
and $\textbf{b}_-$ maps into $W'_{2,-}\otimes_\kappa \OO_{\PP^1}$, the
images are mutually orthogonal.  Therefore $E'_2$ is isotropic.  It is
straightforward to compute that $E'_2$ and $E'_1$ are locally direct
summands of rank $2$, respectively rank $1$.  It is also
straightforward to compute that $R'_3$ is a rank $3$ locally direct
summand that annihilates $E'_1$.  Thus it is all of the annihilator of
$E'_1$.  The cokernels $R'_3/E'_2$ and $E'_2/E'_1$ are invertible
sheaves.  Comparing the degrees of $E'_1$, $E'_2$ and $R'_3$, the
cokernels have degree $0$, thus are isomorphic to $\OO_{\PP^1}$.  Of
course $(E'_1)^\vee$ is isomorphic to the ample invertible sheaf
$\OO_{\PP^1}(2)$.
\end{proof}

\medskip\noindent
If $l=1$, define $W''_0$ to be the zero vector space, and define
$E''_0$ to be the zero sheaf on $\PP^1$.  If $l>1$, 
let $(W''_{4l-4},\beta''_{4l-4})$ and $E''_{2l-2,4l-4}$ be as in
Subsection ~\ref{subsec-app} for $a=l-1$ and $b=2l-2$.  Since $l>1$,
$a$ is positive.  And, of course, $b$ equals $2a$.
Thus
Hypothesis~\ref{hyp-app} holds.  
By Lemma~\ref{lem-app}, $(E''_{2l-2,4l-4})^\vee$ is ample.

\medskip\noindent
Define $(W,\beta)$ to be the orthogonal direct sum of
$(W'_2,\beta'_2)$ and $(W''_{4l-4},\beta''_{4l-4})$, which is just
$(W'_2,\beta'_2)$ if $l$ equals $1$.
Define $R_{2l+1}$ to be the direct sum of
$R'_3$ and $E''_{2l-2,4l-4}$.  Define $E_{2l}$ to
be the direct sum of $E'_2$ and $E''_{2l-2,4l-4}$. 
And define $E_{2l-1}$ to be the direct sum of $E'_1$ and 
$E''_{2l-2,4l-4}$.  

\begin{lem} \label{lem-caseIVa}
\marpar{lem-caseIVa}
Assume $l\geq 2$.  Let $k$ equal $2l$ and let $n$ equal $4l$.  The
flag $E_{2l-1}\subset E_{2l} \subset R_{2l+1}\subset W\otimes_\kappa
\OO_{\PP^1}$ is a 
$(k-1,k,k+1)$-flag parametrized by a morphism $\zeta:\PP^1\rightarrow
M$.  The annihilator of $E_{2l-1}$ equals $R_{2l+1}$.  The cokernels
$R_{k+1}/E_k$ and $E_k/E_{k-1}$ are each isomorphic to $\OO_{\PP^1}$.
And $(E_k/E_{k-1})\otimes E_{k-1}^\vee$ is ample.
\end{lem}

\begin{proof}
This follows by combining Lemma ~\ref{lem-caseIVa'} with
the method of proof of
Lemma~\ref{lem-caseIa}.  The novelty is that $R_{2l+1}$ is not
isotropic.  However, it does equal the annihilator of $E_{2l-1}$ with
respect to $\beta_{2l}$.   
\end{proof}

\begin{prop} \label{prop-caseIVa}
\marpar{prop-caseIVb}
Assume $l\geq 2$.  Let $k$ equal $2l$ and let $n$ equal $4l$.   
The morphism $\zeta:\PP^1 \rightarrow M$ associated to the flag in
Lemma~\ref{lem-caseIVa} is a very twisting family of pointed lines on
$\text{Flag}_k(W,\beta)$.  
\end{prop}

\begin{proof}
This is very similar to the proof of Proposition~\ref{prop-caseIa}.
\end{proof}

\medskip\noindent
\textbf{Case IVb, $k$ odd.}
Let $l\geq 1$ be an integer, let $k$ equal $2l+1$, and let $n$ equal
$2k=4l+2$.  For $(n,k)=(2,1)$, there is no very twisting family of
lines.  Let $(W'_2,\beta'_2)$ be as in ~\ref{subsec-app} for $b=1$.
In particular, $W'_2$ is the direct sum of $W'_{1,+}$ and $W'_{1,-}$.

\medskip\noindent
Let $(W''_{4l},\beta''_{4l})$ and $E''_{2l,4l}$ be as in
Subsection ~\ref{subsec-app} for $a=l$ and $b=2l$.  Since $l\geq
1$, $a$ is positive.  And, of course, $b$ equals $2a$.
Thus
Hypothesis~\ref{hyp-app} holds.  

\medskip\noindent
Define $(W,\beta)$ to be the orthogonal direct sum of
$(W'_2,\beta'_2)$ and $(W''_{4l},\beta''_{4l})$.  Define $R_{2l+2}$ to
be the direct sum of $W'_2$ and $E''_{2l,4l}$.  Define $E_{2l+1}$ to
be the direct sum of $W'_{1,+}$ and $E''_{2l,4l}$.  And define
$E_{2l}$ to be $E''_{2l,4l}$.   

\begin{lem} \label{lem-caseIVb}
\marpar{lem-caseIVb}
Assume $l\geq 1$.  Let $k$ equal $2l+1$ and let $n$ equal $4l+2$.  The
flag $E_{2l}\subset E_{2l+1} \subset R_{2l+2}\subset W\otimes_\kappa
\OO_{\PP^1}$ is a 
$(k-1,k,k+1)$-flag parametrized by a morphism $\zeta:\PP^1\rightarrow
M$.  The annihilator of $E_{2l}$ equals $R_{2l+2}$.  The cokernels
$R_{k+1}/E_k$ and $E_k/E_{k-1}$ are each isomorphic to $\OO_{\PP^1}$.
And $(E_k/E_{k-1})\otimes E_{k-1}^\vee$ is ample.
\end{lem}

\begin{proof}
This is similar to the proof of Lemma~\ref{lem-caseIVa}.
\end{proof}

\begin{prop} \label{prop-caseIVb}
\marpar{prop-caseIVb}
Assume $l\geq 1$.  Let $k$ equal $2l+1$ and let $n$ equal $4l+2$.   
The morphism $\zeta:\PP^1 \rightarrow M$ associated to the flag in
Lemma~\ref{lem-caseIVb} is a very twisting family of pointed lines on
$\text{Flag}_k(W,\beta)$.  
\end{prop}

\begin{proof}
This is similar to the proof of Proposition~\ref{prop-caseIVb}.
\end{proof}

\medskip\noindent
An argument similar to the proof of Claim~\ref{claim-grasslines}
proves the very twisting family can be chosen to be an orbit curve.

\bibliography{my}
\bibliographystyle{amsalpha}

\end{document}